\documentclass[fleqn]{article}

\usepackage[T1]{fontenc}
\usepackage[utf8]{inputenc}
\usepackage{amssymb,amsmath,amsfonts,mathrsfs}
\usepackage{latexsym}
\usepackage[usenames,dvipsnames]{color, xcolor}
\usepackage[figuresright]{rotating}
\usepackage{times}
\usepackage{lmodern}
\usepackage{makeidx}
\usepackage[intoc]{nomencl}
\usepackage{overpic}
\usepackage[font=small]{caption}
\usepackage{subcaption}
\usepackage{textcomp}
\usepackage{nicefrac}
\usepackage{cite}

\usepackage{fancyhdr}  

\usepackage[left=2cm,top=3.5cm,right=2cm,bottom=3.5cm,footskip=45pt, footnotesep=24pt]{geometry}
\usepackage{enumitem}
\setlist{nolistsep} 

\usepackage{tikz}
\usetikzlibrary{external}
\usepackage{pgfplots}
\usepackage{tikz-3dplot} 
\usepackage{tikz-network} 
\pgfplotsset{compat=1.18, colormap name=viridis}
\usetikzlibrary[patterns]
\usetikzlibrary{arrows,positioning,shapes}  
\usetikzlibrary{plotmarks}
\usetikzlibrary{decorations.pathmorphing}
\usetikzlibrary{calc}
\usepgfplotslibrary{colorbrewer}
\usepackage{pgfplots}
\usepgfplotslibrary{groupplots}

\usepackage{amsthm}

\theoremstyle{remark}
\newtheorem{remark}{Remark}

\usepackage[pagebackref=false,breaklinks=true,colorlinks=true,urlcolor=blue,linkcolor=black,citecolor=blue]{hyperref}

\usepackage{empheq}
\usepackage[linesnumbered,ruled,vlined]{algorithm2e}
\SetKwInput{KwInput}{Input}  
\SetKwInput{KwOutput}{Output}
\SetKw{KwBy}{by} 
\DeclareMathOperator*{\argmin}{arg\,min}

\newcommand{\bs}{\boldsymbol}

\newcommand{\uexact}{u}
\newcommand{\uunif}{u_{h}}
\newcommand{\uexactpar}{u^{\sigma}}
\newcommand{\uexactnopar}{u}
\newcommand{\uunifnopar}{u_{h}}

\newcommand{\unopar}{u_{\theta}}
\newcommand{\upar}{u_{\theta}^{\sigma}}

\newcommand{\bsfellt}[2]{\raisebox{0.0\height}{$#1{}\textup{\textsf{\fontfamily{PTSans-TLF}\selectfont \textbf{l}}}$}}
\newcommand{\bsfell}{\mathpalette\bsfellt\relax}

\usepackage{xcolor}
\usepackage{listings}
\lstdefinestyle{terminal}
{
    backgroundcolor=\color{black},
    basicstyle=\scriptsize\color{white}\ttfamily
}

\definecolor{Celeste}{HTML}{0176DE}
\definecolor{Azul}{HTML}{173F8A}
\definecolor{Azuloscuro}{HTML}{03122E}
\definecolor{Amarillo}{HTML}{FEC60D}
\definecolor{Amarillooscuro}{HTML}{E3AE00}
\definecolor{Negro}{HTML}{000000}
\definecolor{GrisA}{HTML}{707070}
\definecolor{GrisB}{HTML}{9C9C9C}
\definecolor{GrisC}{HTML}{C6C6C6}
\definecolor{GrisclaroA}{HTML}{EAEAEA}
\definecolor{GrisclaroB}{HTML}{F0F0F0}
\definecolor{GrisclaroC}{HTML}{F6F6F6}
\definecolor{Blanco}{HTML}{FFFFFF}
\definecolor{Verde}{HTML}{00AA00}
\definecolor{Rojo}{HTML}{F24F4F}

\usepackage{pgfplotstable} 
\usepackage{booktabs}      
\usepackage{siunitx}       
\usepackage{makecell} 
\usepackage{multirow}

\usepackage{tikz}
\usepackage{pgfpages, pgfplots, pgfplotstable}
\usepackage{pgfplots}
\usepackage{tikz-3dplot} 
\usepackage{tikz-network}
\pgfplotsset{compat=1.18, colormap name=viridis}
\usetikzlibrary[patterns]
\usetikzlibrary{arrows,positioning,shapes}  
\usetikzlibrary{plotmarks}
\usetikzlibrary{decorations.pathmorphing}
\usetikzlibrary{calc}
\usepgfplotslibrary{colorbrewer}
\usetikzlibrary{external}
\title{An $r$-adaptive finite element method using neural networks for parametric self-adjoint elliptic problems}

\author{Danilo Aballay$^1$, Federico Fuentes$^2$, Vicente Iligaray$^1$, Ángel J. Omella$^3$,\\ David Pardo$^{3,4,5}$, Manuel A. S\'anchez$^2$, Ignacio Tapia$^1$, Carlos Uriarte$^4$\\[1ex]
\normalsize $^1$Faculty of Engineering, Pontificia Universidad Cat\'olica de Chile, Santiago, Chile\\
\normalsize $^2$Institute for Mathematical and Computational Engineering (IMC),\\ \normalsize Pontificia Universidad Cat\'olica de Chile, Santiago, Chile\\
\normalsize $^3$University of the Basque Country (UPV/EHU), Leioa, Spain\\
\normalsize $^4$Basque Center for Applied Mathematics (BCAM), Bilbao, Spain\\
\normalsize $^5$Basque Foundation for Science (Ikerbasque), Bilbao, Spain
}

\date{\vspace{-3.5ex}}

\begin{document}

\maketitle

\begin{abstract}
This work proposes an $r$-adaptive finite element method (FEM) using neural networks (NNs). The method employs the Ritz energy functional as the loss function, currently limiting its applicability to symmetric and coercive problems, such as those arising from self-adjoint elliptic problems. The objective of the NN optimization is to determine the mesh node locations. For simplicity, these locations are assumed to form a tensor product structure in higher dimensions. The method is designed to solve parametric partial differential equations (PDEs). The resulting parametric $r$-adapted mesh generated by the NN is solved for each PDE parameter instance with a standard FEM.
 Consequently, the proposed approach retains the robustness and reliability guarantees of the FEM for each parameter instance, while the NN optimization adaptively adjusts the mesh node locations.
The construction of FEM matrices and load vectors is implemented such that their derivatives with respect to mesh node locations, required for NN training, can be efficiently computed using automatic differentiation. However, the linear equation solver does not need to be differentiable, enabling the use of efficient, readily available `out-of-the-box' solvers. 
The method's performance is demonstrated on parametric one- and two-dimensional Poisson problems.
\end{abstract}

\section{Introduction}

The finite element method (FEM) is one of the most popular techniques for solving partial differential equations (PDEs). Its strength lies in its rigorous mathematical foundation, which provides a framework for error estimation and convergence analysis. FEM solvers are known for their robustness and reliability~\cite{ciarlet2002finite, hughes2003finite, zienkiewicz2005finite, brenner2008mathematical}. A traditional FEM discretizes the spatial domain into a mesh of elements and constructs an approximated PDE solution using possibly high-order~\cite{babuvska1994p, schwab1996p} piecewise-polynomial functions defined over this mesh. While $h$-adaptivity (refining the mesh size) and $p$-adaptivity (increasing the polynomial order) are well-established in the FEM ~\cite{ainsworth1997posteriori, demkowicz2002fully, rachowicz2006fully, pardo2010multigoal, carstensen2014axioms}, $r$-adaptivity, which involves relocating mesh nodes while preserving the dimensionality and polynomial representability, is a less developed area~\cite{budd2009adaptivity, huang2010adaptive}.
In addition, adaptive FEMs suffer from limitations that become particularly acute in the context of parametric PDEs, especially when singular or highly localized solutions arise naturally. In these cases, generating optimal meshes for a family of problems can be a computationally intensive and time-consuming process.

In recent years, neural networks (NNs) have emerged as useful complementary tools for solving PDEs~\cite{dissanayake1994neural, blechschmidt2021three}. Methods like (Variational) Physics-Informed Neural Networks ((V)PINNs)~\cite{karniadakis2021physics, kharazmi2019variational, kharazmi2021hp, rojas2024robust, uriarte2025optimizing} and Deep Ritz~\cite{yu2018deep, liao2019deep, uriarte2023deep} take advantage of the powerful function approximation capabilities of NNs to represent the PDE solutions (see, e.g., \cite{opschoor2024neural}).

However, NN-based solvers also exhibit important limitations. The training process typically involves solving a highly non-convex optimization problem. Standard optimization algorithms, such as stochastic gradient descent, are prone to getting trapped in local minima, resulting in inaccurate or suboptimal solutions. This lack of guaranteed convergence is a major concern. Furthermore, unlike FEM, there is often a scarcity of rigorous theoretical guarantees regarding the convergence and accuracy of NN-based solvers \cite{goodfellow2016deep, petersen2021topological, uriarte2024solving}. Another significant challenge arises from the need for accurate numerical integration techniques in both (V)PINN and Deep Ritz methods. Whether evaluating the residual at collocation points or computing the energy functional, accurate numerical integration is essential \cite{berrone2022variational, rivera2022quadrature}. This can be computationally demanding, especially in high dimensions, and the choice of integration scheme profoundly impacts the accuracy of the obtained solution and the stability of the training \cite{uriarte2023memory, pal2025solving, taylor2025stochastic}.

This work combines traditional FEM and NN-based solvers to create a more robust and efficient method for solving parametric PDEs. The core idea is to use an NN to perform $r$-adaptivity on the mesh of the FEM used to discretize and solve the problem. This strategy effectively addresses the key limitations of standalone NN or FEM approaches. By employing the well-posed FEM formulation, accurate integration routines, and an efficient solver for the resulting systems of linear equations, we avoid the problem of numerical integration and local minima that hinder the training of NN-based solvers. In particular, our previous work \cite{omella2024r} proposed the use of an NN for $r$-adaptivity via the Ritz energy minimization; however, this work employed an NN to approximate the solution itself, which led to convergence difficulties that hindered the recovery of optimal convergence rates. In contrast, the method proposed here employs the NN solely to adapt the mesh, while the PDE is solved using a classical FEM formulation. This strategy ensures that each iteration yields the global minimum of a discretized quadratic problem that approximates the original PDE solution (see, e.g., \cite{uriarte2025optimizing, baharlouei2025least}), providing both robustness and accuracy. Furthermore, we extend this framework to parametric PDEs by introducing an NN that captures the parametric dependence of the solution on the $r$-adapted mesh. Recent related works propose constructing finite-element-type solutions derived from interpolating NNs \cite{ badia2024finite, badia2025adaptive, badia2025compatible}, while \cite{uriarte2022finite, lee2025finite} propose FEM parametric solvers using NNs.

Herein, our primary contribution is the development of an $r$-adaptive FEM solver for parametric PDEs. The NN represents a parameter-dependent finite-element mesh, where the dimensionality in terms of the PDE coefficients may be large. This dynamic $r$-adaptation concentrates the mesh density on regions where the solution exhibits significant variations, such as sharp gradients, boundary layers, or singularities, leading to substantial improvements in accuracy without requiring a globally fine mesh.

In this work, we consider the Ritz minimization formulation to produce optimal $r$-adapted meshes for parametric PDEs. In particular, this restricts our current formulation to symmetric and coercive problems, such as those arising from self-adjoint elliptic problems. 
Extending our approach to non-symmetric or indefinite problems would likely require handling fine-coarse mesh settings, double-loop strategies (see, e.g., \cite{uriarte2023deep}), or first-order systems of least squares formulations (see, e.g., \cite{bersetche2023deep}), which we leave as future work. Moreover, we restrict our implementation to the lowest-order polynomials ($p=1$) in one- and two-dimensional spatial problems and tensor-product meshes. Extending this implementation to higher-order polynomials ($p>1$) is straightforward, while considering irregular geometries requires a more involved mesh-generation strategy that enables automatic differentiation (AD). 

Our new FEM code, compatible with AD, was implemented using the JAX library~\cite{jax2018github}.
Even though JAX-FEM~\cite{xue2023jax} is an existing JAX-based FEM framework, it is designed for fixed meshes and thus unsuitable for $r$-adaptivity. 
Meanwhile, JAX-SSO \cite{Wu2023}, also based on FEM, is specialized to structural optimization problems involving shells.
In contrast, our implementation treats the mesh node locations---central to the $r$-adaptive strategy---as differentiable parameters. While the assembly process requires AD, we show that the FEM linear solver does not, allowing the use of standard FEM direct or iterative solvers for improved flexibility and efficiency.  

The remainder of the work is as follows. \autoref{section2} describes the Ritz minimization procedure, the considered FEM, and the proposed $r$-adaptive method using NNs. \autoref{section3} discusses our implementation and highlights its key advantages and limitations. \autoref{section4} presents multiple numerical results that demonstrate the effectiveness and accuracy of our method on several benchmark Poisson problems in 1D and 2D. Finally, \autoref{section5} concludes the work and outlines future lines of research.

\section{Mathematical framework}\label{section2}

We first describe the non-parametric case (\autoref{section2.1}) and then extend it to the parametric case (\autoref{section2.2}).

\subsection{Non-parametric case}\label{section2.1}

Let us consider the following abstract variational problem: seek $u\in\mathbb{V}$ such that
\begin{equation}\label{eq:model_abstract_problem}
        b(u,v)=\ell(v),\qquad \forall v\in\mathbb{V},
\end{equation} 
where $\mathbb{V}$ is a real Hilbert space equipped with the norm $\Vert\cdot\Vert_{\mathbb{V}}$,  $b:\mathbb{V}\times\mathbb{V}\longrightarrow\mathbb{R}$ is a symmetric, bounded, and coercive bilinear form, i.e., there exist boundedness and coercivity constants $M\geq\alpha>0$ such that $M \|v\|_{\mathbb{V}}^2\geq b(v,v)\geq\alpha\|v\|_{\mathbb{V}}^2$ for all $v\in\mathbb{V}$, and $\ell:\mathbb{V}\longrightarrow\mathbb{R}$ is a linear and continuous functional. By the Lax-Milgram theorem, \eqref{eq:model_abstract_problem} admits a unique solution $u\in\mathbb{V}$ satisfying $\|u\|_\mathbb{V}\leq\tfrac{1}{\alpha}\|\ell\|_{\mathbb{V}'}$, where $\|\cdot\|_{\mathbb{V}'}$ denotes the norm of the dual space $\mathbb{V'}$.

\begin{remark}
PDEs associated with self-adjoint elliptic differential operators can often be manipulated and formulated weakly in the form of \eqref{eq:model_abstract_problem}, where the underlying bilinear form $b$ will be symmetric (due to the operator being self-adjoint). Linear strictly elliptic problems of this nature will be coercive too \cite{yosida2012functional}. Indeed, the coercivity constant $\alpha$ is also often called the ellipticity constant of the underlying operator. Examples of such problems are Poisson problems being characterized by the Laplace operator, which is self-adjoint and strictly elliptic under certain boundary conditions, i.e., there will be a coercivity constant $\alpha>0$ related to the Poincar\'e constant. 
\end{remark}

\begin{remark}
We note that $\ell$ in the formulation of a PDE in the form \eqref{eq:model_abstract_problem} contains all the `problem data', meaning any forcing terms in the `right-hand side' of the PDE and any boundary conditions (either Dirichlet, Neumann, or Robin boundary conditions). Thus, the estimate  $\|u\|_\mathbb{V}\leq\tfrac{1}{\alpha}\|\ell\|_{\mathbb{V}'}$ translates to a stability estimate reflecting the continuous dependence of the solution on the problem data.
\end{remark}

\subsubsection{Ritz minimization}\label{section2.1.1}

Given $b$ and $\ell$ from \eqref{eq:model_abstract_problem}, define the quadratic \emph{Ritz energy}, $\mathcal{J}:\mathbb{V}\longrightarrow\mathbb{R}$, as 
\begin{equation}
    \mathcal{J}(v) := \frac{1}{2} b(v,v)-\ell(v),\qquad v\in\mathbb{V}.
\end{equation} 
If $u\in\mathbb{V}$ is the unique solution to \eqref{eq:model_abstract_problem}, then $b(u,v)=\ell(v)$ for any $v\in\mathbb{V}$ and, in particular, $b(u,u)=\ell(u)$. 
Thus, by the bilinearity and symmetry of $b$, for \emph{any} $v\in\mathbb{V}$,
\begin{equation}\label{eq:Ritzidentity}
    \begin{aligned}
        \mathcal{J}(v)-\mathcal{J}(u) 
            &=\big(\tfrac{1}{2}b(v,v)-\ell(v)\big)-\big(\tfrac{1}{2}b(u,u)-\ell(u)\big)
            =\tfrac{1}{2}b(v,v)-b(u,v)-\tfrac{1}{2}b(u,u)+b(u,u) \\
            &=\tfrac{1}{2}b(v,v)-\tfrac{1}{2}b(u,v)-\tfrac{1}{2}b(v,u)+\tfrac{1}{2}b(u,u)
            =\tfrac{1}{2}b(v-u,v)-\tfrac{1}{2}b(v-u,u)\\
            &=\tfrac{1}{2}b(v-u,v-u)=\tfrac{1}{2} \Vert v-u\Vert^2_b\geq0\,,
    \end{aligned}
\end{equation}
where $\Vert\cdot\Vert_b^2:=b(\cdot,\cdot)$ is the norm in $\mathbb{V}$ induced by the bilinear form (norm-equivalent to $\|\cdot\|_{\mathbb{V}}$ due to the coercivity and boundedness of $b$). 
Thus, $\mathcal{J}(v)>\mathcal{J}(u)$ for any $v\neq u$, meaning the unique solution $u\in\mathbb{V}$ of \eqref{eq:model_abstract_problem} is also the unique minimizer of $\mathcal{J}$ over $\mathbb{V}$, i.e., $u=\argmin_{v\in\mathbb{V}}\mathcal{J}(v)$.

\subsubsection{Discretization and \texorpdfstring{$r$-adaptivity}{r-adaptivity}}
\label{section2.1.2}

For a 1D spatial domain $\Omega=(a,b)$, we consider a set of trainable variables $\theta=(\theta_1,\theta_2,\ldots,\theta_n) \in\mathbb{R}^{n}$ and apply the \textit{softmax} activation function to obtain a partition of unity $\delta (\theta)=(\delta_1,\delta_2,\ldots,\delta_{n})$ as:
\begin{subequations} \label{eq:1Dmeshconstruction}
\begin{align}
    \delta_i &=\displaystyle \frac{\exp({\theta_i})}{\sum_{j=1}^{n}\exp(\theta_j)}\in (0,1), \qquad 1\leq i\leq n,\\
    1 &= \sum_{i=1}^n \delta_i.
\end{align} 
This partition serves to define the nodal points in $\Omega$ as follows:
\begin{alignat}{2}
    x_0 &:=a,\label{eq:a} &&\\
    x_i &:= x_{i-1} + (b-a) \textstyle  \delta_i,\qquad &&1\leq i\leq n.\label{eq:xi}
\end{alignat} Note that by construction $x_n = b$.

In addition to the above parameterized nodal points for adaptivity, we may consider some additional fixed (non-adaptive) interior nodes $a<x_{n+1}<x_{n+2}<\ldots<x_N<b$ (used, for example, to define material interfaces). Hence, the final sorted list of adaptive and non-adaptive nodal points defines our FEM mesh with $N$ elements,
\begin{alignat}{2}\label{eq:nodesorting}
    \bs{\mathsf{x}}_\theta &:= (\tilde{x}_0,\ldots,\tilde{x}_N) 
        =\texttt{sort}\Big[x_0,\underbrace{x_1,x_2,\ldots,x_{n-1}}_{\text{adaptive nodes}},x_n,\underbrace{x_{n+1},x_{n+2},\ldots,x_N}_{\text{fixed nodes}}\Big].
\end{alignat} 
\end{subequations}  

We consider piecewise-linear FEM `hat functions' $\bs{\mathsf{v}}_{\theta}=\{v_0, \ldots, v_N\}$ associated with $\bs{\mathsf{x}}_\theta$ as shown in \autoref{fig:tent_functions}, satisfying $v_i(\tilde{x}_j)=\delta_{ij}$ where $\delta_{ij}$ is a Kronecker delta. One can consider a much larger set of FEM functions (e.g., piecewise polynomials of higher order), but for simplicity we limit ourselves to the most straightforward case. 
\begin{figure}[htbp]
    \centering
\begin{tikzpicture} 
        \draw[dotted, thick][-] (0,0) -- (10.5,0) ; 
		\draw (0,-0.35) node{$a=\tilde{x}_0$};
		\draw[color=red, thick][-]  (1.5,0) -- (0,1.5);
		\draw (0,1.9) node{\color{red}$v_{0}(x)$};
		\draw (1.5,-0.35) node{$\tilde{x}_1$};
		\draw (3,-0.35) node{$\tilde{x}_2$};
        \draw (4.5,-0.35) node{$\tilde{x}_3$};
        \draw (6,-0.35) node{$\tilde{x}_4$};
        \draw (7.5,-0.35) node{$\ldots$};
        \draw (9,-0.35) node{$\tilde{x}_{N-1}$};
		\draw (10.5,-0.35) node{$\tilde{x}_{N}=b$};
		\draw[thick][-]  (0,0) -- (1.5,1.5) node[anchor=north] {};
		\draw[thick][-]  (1.5,1.5) -- (3,0) node[anchor=north] {};
		\draw (1.5,1.9) node{$v_{1}(x)$};
		\draw[color=green!40!black, thick][-]  (1.5,0) -- (3,1.5) node[anchor=north] {};
		\draw[color=green!40!black, thick][-]  (3,1.5) -- (4.5,0) node[anchor=north] {};
		\draw (3,1.9) node{\color{green!40!black}$v_2(x)$};		
		\draw[color=blue, thick][-]  (3,0) -- (4.5,1.5) node[anchor=north] {};
		\draw[color=blue, thick][-]  (4.5,1.5) -- (6,0) node[anchor=north] {};
		\draw (4.5,1.9) node{\color{blue}$v_3(x)$};
		\draw[color=purple, thick][-]  (9,0) -- (10.5,1.5);
		\draw (10.5,1.9) node{\color{purple}$v_{N}(x)$}; node[anchor=north] {};
	\foreach \x in {0,1.5,3,4.5,6,9,10.5}{
         \node[circle, draw=black, fill=black,anchor=north, inner sep=1.5pt, minimum size=0.5pt] (a) at (\x, 0.08){};
         }
\end{tikzpicture}
    \caption{Piecewise-linear FEM basis functions in $\Omega=(a,b)$.}
    \label{fig:tent_functions}
\end{figure}

Next, from the boundary conditions, we construct a discrete FEM space as the span of the relevant subset of $\{v_0, \ldots, v_N\}$, which we denote by $\mathbb{V}_\theta$ to reflect its dependence on the trainable variables $\theta=(\theta_1,\theta_2,\ldots,\theta_n)$.
More precisely, we denote $\Lambda_{\mathrm{dof}}\subseteq\{0,\ldots,N\}$ as the indices associated with the nodes \emph{without} assigned Dirichlet boundary conditions, and $\Lambda_{D}=\{0,\ldots,N\}\smallsetminus\Lambda_{\mathrm{dof}}$ as those indices where the corresponding nodes have assigned Dirichlet boundary conditions. 
Thus, $\mathbb{V}_\theta=\mathrm{span}\big(\{v_j\mid j\in\Lambda_{\mathrm{dof}}\}\big)$ and the FEM solution $u_\theta\in\mathbb{V}_\theta$ to \eqref{eq:model_abstract_problem} may be written as
\begin{equation}\label{eq:nonparametric-solution}
   u_\theta(x) = \sum_{j\in\Lambda_{\mathrm{dof}}} c_j \; v_j(x)\,,
\end{equation}
where $\bs{\mathsf{c}}_\theta=(c_j)_{j\in\Lambda_{\mathrm{dof}}}$ is a vector of coefficients associated with the basis $\{v_j\}_{j\in\Lambda_{\mathrm{dof}}}$, and it is given by $\bs{\mathsf{c}}_\theta=\bs{\mathsf{B}}^{-1}_\theta\bsfell_\theta$ with 
$(\bs{\mathsf{B}}_{\theta})_{ij}:=b(v_{j},v_{i})$ and $(\bsfell_{\theta})_j:=\ell(v_j)$ for $i,j\in\Lambda_{\mathrm{dof}}$.
There is some abuse of notation in that $(\bs{\mathsf{B}}_{\theta})_{ij}$ does not exist when $i$ or $j$ is in $\Lambda_{D}$ (so $\bs{\mathsf{B}}_{\theta}$ is \emph{not} an $(N+1)\times (N+1)$ matrix), but it should be clear that the $|\Lambda_{\mathrm{dof}}|\times|\Lambda_{\mathrm{dof}}|$ matrix $\bs{\mathsf{B}}_{\theta}$ and the $|\Lambda_{\mathrm{dof}}|\times1$ vector $\bsfell_{\theta}$ are defined only by the indices in $\Lambda_{\mathrm{dof}}$ (which technically should require a proper relabeling to indices in $\{1,\ldots,|\Lambda_{\mathrm{dof}}|\}$, but we omit these details for the sake of brevity). 
We also note that the Dirichlet boundary conditions enter the right-hand side $\ell(v)$ in \eqref{eq:model_abstract_problem}, which will contain the term $-b(u_0,v)$, where $u_0(x)=\sum_{j\in\Lambda_D} c_j \; v_j(x)$ has the $c_j$ selected to satisfy the desired Dirichlet boundary conditions, and the solution of interest (i.e., the one that solves a PDE) will then take the form $\tilde{u}(x)=u_0(x)+u_\theta(x)=\smash{\sum_{j=0}^N} c_jv_j(x)$.

For a 2D rectangular domain $\Omega=(a^x,b^x)\times(a^y\times b^y)$, we restrict to tensor products of 1D meshes and proceed analogously, so that a set of trainable variables is given by $\theta=(\theta^x,\theta^y)\in\mathbb{R}^{n_x+n_y}$, where $\theta^x=(\theta_1^x,\theta_2^x,\ldots,\theta_{n_x}^x)$ and $\theta^y=(\theta_1^y,\theta_2^y,\ldots,\theta_{n_y}^y)$.
Then, we apply the softmax function to obtain $\delta^x(\theta^x)$ and $\delta^y(\theta^y)$ and proceed as in \eqref{eq:1Dmeshconstruction} to obtain a set of $N_xN_y$ quadrilateral elements characterized by $(N_x+1)(N_y+1)$ 2D nodes $\bs{\mathsf{x}}_\theta=\big((\tilde{x}_0,\tilde{y}_0),\ldots,(\tilde{x}_{N_x},\tilde{y}_{N_y})\big)$. With these nodes, we can then construct the associated basis of $(N_x+1)(N_y+1)$ bilinear quadrilateral hat functions $\bs{\mathsf{v}}_{\theta}$, define the sets $\Lambda_{\mathrm{dof}}$ and $\Lambda_{D}$ according to the boundary conditions, and proceed as described above. 

Given an initial set of weights $\theta^{(0)}$, we iterate for each $t=0,1,\ldots,T$ according to the following four-step $r$-adaptive procedure:
\begin{enumerate}[label=(\roman*),itemsep=1mm]
    \item Given $\theta^{(t)}$, generate the mesh $\bs{\mathsf{x}}^{(t)}_\theta$.
    \item Construct $\bs{\mathsf{B}}^{(t)}_\theta$ and $\bsfell^{(t)}_\theta$.
    \item Solve $\bs{\mathsf{B}}^{(t)}_\theta \bs{\mathsf{c}}^{(t)}_\theta = \bsfell^{(t)}_\theta$.
    \item Update $\theta^{(t)}$ using a gradient-descent-based algorithm applied to the Ritz loss function,
    \begin{equation}
        \theta^{(t+1)} = \theta^{(t)} - \eta^{(t)}\nabla_{\theta}\mathcal{J}(u_\theta^{(t)}).
    \end{equation}
    Here, $\eta^{(t)}$ denotes the iteration-dependent learning rate, which may come from sophisticated gradient-descent-based optimizers like Adam \cite{kingma2014adam}, and $\nabla_{\theta}\mathcal{J}(u_\theta^{(t)})$ is the gradient with respect to $\theta$ of the Ritz functional evaluated as
    \begin{equation}\label{eq:loss_with_theta}
        \mathcal{J}(u_\theta^{(t)}) = \frac{1}{2}\bs{\mathsf{B}}^{(t)}_\theta\bs{\mathsf{c}}^{(t)}_\theta\cdot\bs{\mathsf{c}}^{(t)}_\theta -\bsfell^{(t)}_\theta\cdot\bs{\mathsf{c}}^{(t)}_\theta\,.
    \end{equation}
\end{enumerate} 
\autoref{fig:diagram_adaptivity_nonparametric} shows a flowchart  of the described $r$-adaptive method.

\begin{figure}[htbp]
    \centering
    \resizebox{0.9\textwidth}{!}{%
    \input{figures/tikz/diagram_adativity_nonparametric}%
    }
    \caption{Flowchart of the proposed $r$-adaptive method in the non-parametric case.}
    \label{fig:diagram_adaptivity_nonparametric}
\end{figure}

\begin{remark}[No need to compute derivatives of the solver of linear equations] \label{remark:vanishing_partial}
Viewing $\mathcal{J}=\mathcal{J}(\bs{\mathsf{B}}_\theta,\bsfell_\theta,\bs{\mathsf{c}}_\theta)$, notice
\begin{equation}
    \nabla_\theta\mathcal{J}
        =\frac{\partial\mathcal{J}}{\partial\bs{\mathsf{B}}_\theta}
        \frac{\partial\bs{\mathsf{B}}_\theta}{\partial\theta}
            +\frac{\partial\mathcal{J}}{\partial\bsfell_\theta}
            \frac{\partial\bsfell_\theta}{\partial\theta}
                +\frac{\partial\mathcal{J}}{\partial\bs{\mathsf{c}}_\theta}
                \frac{\partial\bs{\mathsf{c}}_\theta}{\partial\theta}
        =\frac{1}{2}\bs{\mathsf{c}}_\theta^\mathsf{T}\frac{\partial\bs{\mathsf{B}}_\theta}{\partial\theta}\bs{\mathsf{c}}_\theta
            -\frac{\partial\bsfell_\theta}{\partial\theta}\cdot\bs{\mathsf{c}}_\theta\,,
\end{equation}
since $\partial\mathcal{J}/\partial{\bs{\mathsf{c}}_\theta} = \bs{\mathsf{B}}_\theta\, \bs{\mathsf{c}}_\theta-\bsfell_\theta = 0$. Therefore, the derivative $\partial\bs{\mathsf{c}}_\theta/\partial\theta$ need not be computed, and one only needs access to $\bs{\mathsf{c}}_\theta = \bs{\mathsf{B}}_\theta^{-1}\bsfell_\theta$, which could be solved outside any AD framework using a preferred solver of linear equations for the FEM. 
This decoupling allows the integration of standard, potentially non-differentiable, high-performance FEM solvers to be included within the overall optimization loop.
\end{remark}

\begin{remark}\label{remark:alternative_expressions}
    Instead of $\mathcal{J}(u_\theta)=\frac{1}{2}\bs{\mathsf{B}}_\theta\bs{\mathsf{c}}_\theta\cdot\bs{\mathsf{c}}_\theta -\bsfell_\theta\cdot\bs{\mathsf{c}}_\theta$ in \eqref{eq:loss_with_theta}, alternative reformulations, such as $\mathcal{J}(u_\theta)=-\tfrac{1}{2}\bs{\mathsf{B}}_\theta\bs{\mathsf{c}}_\theta\cdot\bs{\mathsf{c}}_\theta$ or $\mathcal{J}(u_\theta)=-\tfrac{1}{2}\bsfell_\theta\cdot\bs{\mathsf{c}}_\theta$, could be used. However, these forms are less advantageous because their partial derivative with respect to $\bs{\mathsf{c}}_\theta$ generally does not vanish (i.e., $\partial\mathcal{J}/\partial\bs{\mathsf{c}}_\theta \neq 0$). Consequently,  \autoref{remark:vanishing_partial} is not applicable. Thus, whenever possible, we recommend using $\mathcal{J}(u_\theta)$ in \eqref{eq:loss_with_theta} as the objective.
\end{remark}

\subsection{Parametric case}\label{section2.2}

For each $\sigma$ belonging to a parameter space $\Sigma$, let us consider the following variational problem: seek $u^\sigma\in\mathbb{V}$ such that 
\begin{equation}\label{eq:model_abstract_problem_parametric}
    b^\sigma(u^\sigma,v)=\ell^\sigma(v),\qquad \forall v\in\mathbb{V},
\end{equation} 
where $\mathbb{V}$ is a real Hilbert space equipped with the norm $\Vert\cdot\Vert_{\mathbb{V}}$, $b^\sigma:\mathbb{V}\times\mathbb{V}\longrightarrow\mathbb{R}$ is a symmetric, bounded, and coercive $\sigma$-dependent bilinear form, and $\ell^\sigma:\mathbb{V}\longrightarrow\mathbb{R}$ is a linear and continuous $\sigma$-dependent functional. Then, each $\sigma$-dependent problem in \eqref{eq:model_abstract_problem_parametric} admits a unique solution $u^\sigma=\argmin_{v\in\mathbb{V}}\mathcal{J}^\sigma(v)$, where the \emph{$\sigma$-dependent Ritz energy} is defined as
\begin{equation}
    \mathcal{J}^\sigma(v) := \frac{1}{2} b^\sigma(v,v)-\ell^\sigma(v),\qquad v\in\mathbb{V}\,.
\end{equation} We emphasize that the range of $\mathcal{J}^\sigma$ can vary widely as a function of $\sigma$. Therefore, one might desire to balance the above so as to produce qualitatively equivalent $\sigma$-dependent functionals whose minima are close or even coincide. In a first approach, one might propose 
\begin{equation}
    \frac{\mathcal{J}^\sigma(v)}{|\mathcal{J}^\sigma(u^\sigma)|} = \frac{\frac{1}{2}b^\sigma(v,v)-\ell^\sigma(v)}{\vert \min_{w\in\mathbb{V}} \frac{1}{2}b^\sigma(w,w)-\ell^\sigma(w)\vert},\qquad v\in\mathbb{V},
\end{equation} 
because $\mathcal{J}^\sigma(u^\sigma)/|\mathcal{J}^\sigma(u^\sigma)|=\min_{v\in\mathbb{V}} \mathcal{J}^\sigma(v)/|\mathcal{J}^\sigma(u^\sigma)| = -1$ for any $\sigma\in\Sigma$ whenever $\mathcal{J}^{\sigma}(u^{\sigma})\neq 0$, which occurs when $u^{\sigma}\neq 0$. However, this requires knowing the value of $\mathcal{J}^\sigma(u^\sigma)$ beforehand, which is often precisely what one wishes to compute with precision. With this in mind, we can instead use a proxy that is cheaply computed and that aims to achieve the same goal of balancing $\mathcal{J}^\sigma$ on different values of $\sigma$. For $N\in \mathbb{N}$ and $h=\tfrac{1}{N}>0$, we propose considering the \emph{$\sigma$-balanced Ritz energy} as 
\begin{equation}
    \tilde{\mathcal{J}}^\sigma(v) := \frac{\mathcal{J}^\sigma(v)}{|\mathcal{J}^\sigma(u_h^\sigma)|},\qquad v\in\mathbb{V}\,,
    \label{eq:balancedRitzdef}
\end{equation}
where $u_h^\sigma=\argmin_{v_h\in\mathbb{V}_h}\mathcal{J}^\sigma(v_h)$ is the `discrete' solution to the minimization problem over a finite-dimensional subspace $\mathbb{V}_h\subseteq\mathbb{V}$ coming from a \textit{uniform} discretization of size $h=\tfrac{1}{N}$ of the spatial domain $\Omega$. As a result, we expect $|\tilde{\mathcal{J}}^{\sigma}(u^{\sigma})|=\mathcal{O}(1)$ across different values of $\sigma$.

\begin{remark}\label{rmk:negativeRitz}
    From \autoref{remark:alternative_expressions}, we have that $\mathcal{J}^\sigma(u_h^{\sigma}) = -\tfrac{1}{2}(\bs{\mathsf{c}}_h^\sigma)^\mathsf{T}\bs{\mathsf{B}}_h^\sigma\bs{\mathsf{c}}_h^\sigma=-\tfrac{1}{2}b^\sigma(u_h^{\sigma},u_h^{\sigma})<0$ as long as $u_h^{\sigma}\neq0$. Therefore, $\mathcal{J}^\sigma(u_h^{\sigma})$ is `always' negative and $\tilde{\mathcal{J}}^\sigma$ is well defined.
\end{remark}

\subsubsection{Discretization and \texorpdfstring{$r$-adaptivity}{r-adaptivity}}
\label{section2.2.2}

For a 1D spatial domain $\Omega=(a,b)$, we consider the nodal-point discretization $\bs{\mathsf{x}}_\theta^\sigma$ to be $\sigma$-dependent by first applying a fully-connected NN to $\sigma\in\Sigma$ and then utilizing our previous mesh generation scheme:
\begin{subequations}
\begin{alignat}{2}
    \mathsf{z}_0 &:= \sigma\in\Sigma, && \label{eq:net1}\\
    \mathsf{z}_l(\sigma) &:= \varphi(\mathsf{W}_l \; \mathsf{z}_{l-1}(\sigma) + \mathsf{b}_l),\qquad && 1\leq l\leq L-1,\label{eq:net2}\\
    \delta(\sigma) &:= (\delta_1(\sigma),\delta_2(\sigma),\ldots,\delta_n(\sigma)) = \text{softmax}(\mathsf{W}_L\; \mathsf{z}_L(\sigma)), &&\label{eq:softmax}\\
    x_0 &:=a, &&\\
    x_i(\sigma) &:= x_{i-1}(\sigma) + (b-a)\delta_i(\sigma),\qquad &&1\leq i\leq n,\\
    \bs{\mathsf{x}}_\theta^\sigma &:= \texttt{sort}\left[x_0,x_1(\sigma),x_2(\sigma),\ldots,x_{n-1}(\sigma),x_n,x_{n+1},\ldots,x_N\right]\label{eq:xi_sigma}. \qquad &&
\end{alignat}
\end{subequations} 
Here, $\theta=(\mathsf{W}_1, \mathsf{b}_1, \mathsf{W}_2, \mathsf{b}_2,\ldots, \mathsf{W}_{L-1}, \mathsf{b}_{L-1}, \mathsf{W}_L)$ is the set of trainable variables (weights and biases) of the NN, $L$ is its depth, and $\varphi$ is an activation function that acts componentwise.

As a result, $\theta$ delivers $\sigma$-dependent meshes $\bs{\mathsf{x}}_\theta^\sigma$, which we can subordinate to FEM bases $\bs{\mathsf{v}}_{\theta}^\sigma=\{v_0^\sigma, \ldots, v_N^\sigma\}$, as described in \autoref{section2.1.2}.
The $\sigma$-dependent discrete finite dimensional spaces $\mathbb{V}_\theta^\sigma=\mathrm{span}\big(\{v_j^\sigma : j\in\Lambda_{\mathrm{dof}}^\sigma\}\big)$ are then defined after having properly labeled the nodes associated to Dirichlet boundary conditions in $\Lambda_D^\sigma$, leaving the remaining indices as $\Lambda_{\mathrm{dof}}^\sigma$.
The optimal FEM solution to \eqref{eq:model_abstract_problem_parametric} in $\mathbb{V}_\theta^\sigma$ is then given by 
\begin{equation}
    u_\theta^\sigma(x)=\sum_{j\in\Lambda_{\mathrm{dof}}^\sigma} c_j^\sigma\,v_j^\sigma(x)\,,
\end{equation}
where $\bs{\mathsf{c}}_\theta^\sigma=(c_j^\sigma)_{j\in\Lambda_{\mathrm{dof}}}$ is a vector of coefficients associated with the basis $\{v_j^\sigma\}_{j\in\Lambda_{\mathrm{dof}}^\sigma}$, and is given by $\bs{\mathsf{c}}_\theta^\sigma=(\bs{\mathsf{B}}_\theta^\sigma)^{-1}\bsfell_\theta^\sigma$ with 
$(\bs{\mathsf{B}}_{\theta}^\sigma)_{ij}:=b^\sigma(v_{j}^\sigma,v_{i}^\sigma)$ and $(\bsfell_{\theta}^\sigma)_j:=\ell^\sigma(v_j^\sigma)$ for $i,j\in\Lambda_{\mathrm{dof}}^\sigma$.
The generalization to 2D and 3D tensor-product meshes is straightforward (see \autoref{section2.1.2}).

The space $\mathbb{V}_\theta^\sigma$ is conforming to a larger infinite-dimensional space $\mathbb{V}$, so $\mathbb{V}_\theta^\sigma\subseteq\mathbb{V}$, and the solution $u^\sigma$ to \eqref{eq:model_abstract_problem_parametric} satisfies the following optimality condition in relation to its $\mathbb{V}_\theta^\sigma$-restricted solution $u_\theta^\sigma\in\mathbb{V}_\theta^\sigma$, 
\begin{equation}\label{eq:inequality}
    \tilde{\mathcal{J}}^\sigma(u^\sigma) = \min_{v\in\mathbb{V}} \tilde{\mathcal{J}}^\sigma(v) \leq  \min_{v\in\mathbb{V}^\sigma_\theta} \tilde{\mathcal{J}}^\sigma(v) = \tilde{\mathcal{J}}^\sigma(u_\theta^\sigma),\qquad \text{for every } \sigma\in\Sigma \text{ and any }\theta.
\end{equation} 
Then, we parameterize the balanced Ritz functional across all $\sigma$ as follows:
\begin{equation}\label{eq:optimal_Ritz_global_parameterized}
    \theta\mapsto\tilde{\mathcal{J}}^\Sigma(\theta) := \int_\Sigma \tilde{\mathcal{J}}^\sigma(u^\sigma_\theta)\; d\mu(\sigma),\qquad \tilde{\mathcal{J}}^\sigma(\;\cdot\;)=\frac{\mathcal{J}^\sigma(\;\cdot\;)}{|\mathcal{J}^\sigma(u_h^\sigma)|},
\end{equation} 
where $\mu$ is an appropriate measure for $\Sigma$, and $u_h^\sigma=\argmin_{v_h\in\mathbb{V}_h}\mathcal{J}^\sigma(v_h)$ is the solution at the uniform mesh of size $h=\tfrac{1}{N}$ (with $h$ selected to be compatible with fixed material interfaces, when present).
By \eqref{eq:inequality}, it is immediate that $\tilde{\mathcal{J}}_{\mathbb{V}}^\Sigma:=\int_\Sigma \tilde{\mathcal{J}}^\sigma(u^\sigma)\; d\mu(\sigma)$ is a lower bound for $\tilde{\mathcal{J}}^\Sigma(\theta)$ in \eqref{eq:optimal_Ritz_global_parameterized}, i.e., $\tilde{\mathcal{J}}_{\mathbb{V}}^\Sigma \leq \tilde{\mathcal{J}}^\Sigma(\theta)$ for any $\theta$.

In practice, for each iteration of the gradient-descent optimization, we approximate the above integral via the following finite weighted average:
\begin{equation}\label{eq:loss_parametric}
\tilde{\mathcal{J}}^\Sigma(\theta) \approx \frac{1}{J} \sum_{j=1}^J \tilde{\mathcal{J}}^{\sigma_j}(u_\theta^{\sigma_j}),
\end{equation} 
where $\{\sigma_1,\sigma_2,\ldots,\sigma_J\}\subset\Sigma$ is a finite sample of parameters coming from some problem-dependent distribution (in the simplest case, a uniform one).

\autoref{fig:diagram_adaptivity_parametric} shows the flowchart for our $r$-adaptive method for parametric problems. 

\begin{figure}[htbp]
    \centering
    \resizebox{0.98\textwidth}{!}{%
    \input{figures/tikz/diagram_adaptivity_parametric}%
    }
    \caption{Flowchart of the proposed $r$-adaptive method in the parametric case.}
    \label{fig:diagram_adaptivity_parametric}
\end{figure}

\subsection{On exact computation of integrals within \texorpdfstring{$r$-adaptivity}{r-adaptivity}}
\label{sec:quadraturevsexactintegration}

Before proceeding further, we alert the reader on a delicate issue regarding the accurate computation of the Ritz energy functional $\mathcal{J}(u)$, which translates to evaluating the forms $b(u,u)$ and $\ell(u)$, normally containing integrals.
In FEM, these are usually calculated elementwise via standard quadrature rules, such as Gauss-Legendre, which are exact for typical polynomial FEM integrands, but only \emph{approximate} non-polynomial ones which may appear in $b(u,u)$ or $\ell(u)$ under certain material heterogenities and forcings.
While this quadrature error diminishes with $h$-refinements \cite{brenner2008mathematical}, it is not negligible for $r$-adaptivity, as we show next.

\begin{figure}[htbp]
    \centering
    \includegraphics[width=0.98\linewidth]{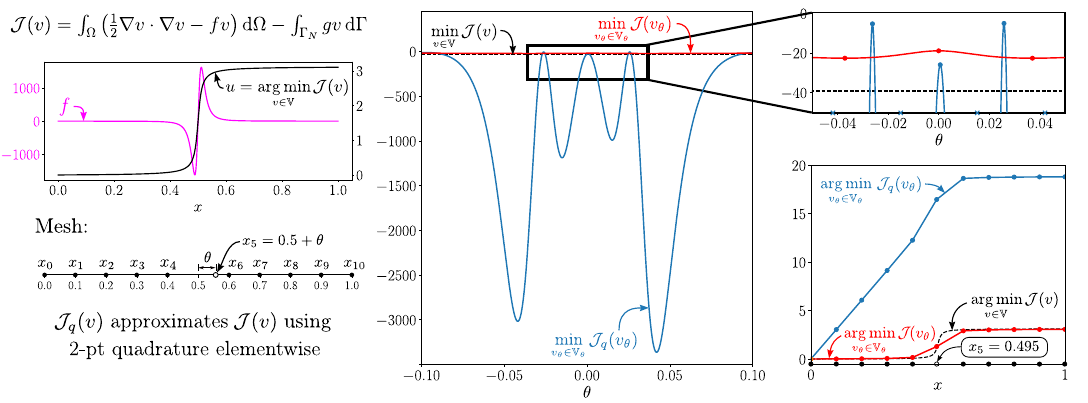}
    \caption{
        For a Poisson problem $-u''(x)=f(x)$ with solution $u(x)=\arctan(\alpha(x-s))+\arctan(\alpha s)$ for $\alpha=50$ and $s=0.5$, as shown in the left panel, we illustrate the landscapes of the minimum Ritz energy over piecewise-linear functions as a single mesh-node position is varied. The minimal Ritz energy landscapes corresponding to the exact and inexact integration and their associated minimizers are shown in red and blue, respectively, in the central and right panels.
    }
    \label{fig:QuadratureVsAnalytic}
\end{figure}

Consider solving $-u''(x)=f(x):=\tfrac{2\alpha^3(x-s)}{(1+\alpha^2(x-s)^2)^2}$ in $\Omega=(0,1)$ with $u(0)=0$ and $u'(1)=g:=\tfrac{\alpha}{1+\alpha^2(1-s)^2}$, whose exact solution is $u(x)=\arctan(\alpha(x-s))+\arctan(\alpha s)$.
As noted in \autoref{section2.1.1}, this solution is the unique minimizer of the Ritz functional $\mathcal{J}(v)=\tfrac{1}{2}b(v,v)-\ell(v)$, where $b(v,v):=\int_\Omega|\nabla v|^2\,\mathrm{d}\Omega$ and $\ell(v):=\int_\Omega f v\,\mathrm{d}\Omega-v(1)g$ for $v\in\mathbb{V}:=\{v\in H^1(\Omega) : v(0)=0\}$.
Now, consider mesh nodes $0=x_0<x_1<\ldots<x_N=1$, elements $\Omega^e=(x_{e-1},x_e)$ for $e=1,\ldots,N$, and corresponding discrete space $\mathbb{V}_\theta=\mathrm{span}\big(\{v_1,\ldots,v_N\}\big)\subseteq\mathbb{V}$ spanned by the associated hat functions $v_j$, as described in \autoref{section2.1.2}.
For simplicity, and to illustrate the main point, let $\alpha=50$, $s=0.5$, and $N=10$, and fix all nodal points $x_j$ at a uniform mesh (i.e., $x_j=\tfrac{j}{10}$) except $x_5=0.5+\theta$, which will then be a $\theta$-dependent $r$-adaptable node that modifies $\mathbb{V}_\theta$.
Thus, $\min_{v_\theta\in\mathbb{V}_\theta}\mathcal{J}(v_\theta)$ will vary as a function of $\theta$, whose landscape is depicted in red in \autoref{fig:QuadratureVsAnalytic}.
This curve, showing two local minima in the top-right panel, lies entirely above the true global minimum $\mathcal{J}(u)$.
An $r$-adaptive strategy optimizing $\theta$ should approach one of these local minima. However, this relies on being able to compute $\mathcal{J}(v_\theta)$ exactly.

Let $\mathcal{J}_q(v_\theta)$ be the \emph{approximation} of $\mathcal{J}(v_\theta)$ resulting from computing the elementwise integrals in $\mathcal{J}(v_\theta)$ for $v_\theta\in\mathbb{V}_\theta$ with two-point Gauss-Legendre quadrature.
Clearly, $\mathcal{J}_q(v_\theta)=\tfrac{1}{2}b(v_\theta,v_\theta)-\ell_q(v_\theta)$ since $b(v_\theta,v_\theta)$ is computed exactly.
However, $\ell_q(v_\theta)\neq\ell(v_\theta)$ and this has a dramatic effect on the Ritz energy landscape, which is illustrated in blue in \autoref{fig:QuadratureVsAnalytic}.
This new landscape has different local minima, and, notably, those minima lie orders of magnitude below the theoretical global minimum $\mathcal{J}(u)=\min_{v\in \mathbb{V}}\mathcal{J}(v)$.
Note that a small change in $\theta$, for instance, $\theta=-0.005$, results in a huge jump in $\min_{v_\theta\in\mathbb{V}_\theta}\mathcal{J}_q(v_\theta)$, when compared to $\min_{v_\theta\in\mathbb{V}_\theta}\mathcal{J}(v_\theta)$.
If a NN is used to solve this problem, the underlying optimizer will quickly approach the minimizers in the degraded landscape, instead of the theoretical ones. Moreover, the resulting minimizer will be quite far from the theoretical optimizer, as illustrated in the bottom-right panel of \autoref{fig:QuadratureVsAnalytic}. 

Therefore, it is crucial to evaluate $\mathcal{J}(v_\theta)$ accurately.
One solution is to use highly-accurate quadrature rules and hope that the resulting $\mathcal{J}_q(v_\theta)$ is sufficiently close to $\mathcal{J}(v_\theta)$. 
Alternatively, if possible, one may implement analytical integration routines to evaluate $\mathcal{J}(v_\theta)$ exactly, or limit oneself to problems with suitable piecewise-polynomial integrands that can be evaluated exactly with (sufficiently accurate) numerical quadrature.
In this work, most non-polynomial expressions for the integrands lent themselves to analytic integration, so, when possible, we implemented specialized routines to integrate them exactly and avoid this potential issue.

\section{Implementation of JAX-based finite element code}
\label{section3}

\subsection{Automatic differentiation and just-in-time compilation}

The code was implemented in Python using the library JAX \cite{jax2018github,frostig2018compiling} to: (a) allow for automatic differentiation (AD) of its components, which is necessary in order to train the NN using stochastic gradient descent, and (b) enable XLA (Accelerated Linear Algebra) usage via just-in-time (JIT) compilation. 
It is available at the repository rFEM-NN \cite{rFEMgithub}, which contains all the details described below and in the numerical experiments of \autoref{section4}.

As mentioned in the introduction, there are other JAX-based FEM codes such as JAX-FEM \cite{xue2023jax} and JAX-SSO \cite{Wu2023}. JAX-SSO \cite{Wu2023} is specialized to shape and topology optimization problems involving shells. Meanwhile, JAX-FEM \cite{xue2023jax} is developed with inverse problems in mind, often involving nonlinear forward problems coming from computational mechanics, and is built to solve problems efficiently on a fixed given mesh (using JAX's AD). Our $r$-adaptive method, however, requires differentiating the loss function with respect to the parameters that determine the mesh node positions. Because this implies optimizing a variable mesh geometry, it fundamentally changes how FEM basis functions must be handled, making it incompatible with JAX-FEM's fixed-mesh structure. This forced us to build our implementation from the ground up.

To benefit from JAX's JIT compilation, a careful implementation is required.
Indeed, control flow constructs such as \texttt{if} statements or loops can interfere with tracing when their conditions or iteration bounds depend on runtime values (e.g., JAX arrays), resulting in errors or recompilation. 
Thus, most \texttt{if} statements disallow exploiting the advantages of JIT compilation. To handle these situations, one must use a suite of JAX functions that are compatible with both JIT and AD.
This presents some challenges, but despite this limitation, the existing suite is sufficient to achieve the desired objectives. Further implementation details follow.

\subsection{Dirichlet node labeling}
\label{sec:Meshgenerationandlabeling}

The mesh is composed of both fixed nodes (e.g., those needed for interfaces) and NN-trained adaptable nodes (see \autoref{section2.1.2}). After each update to the adaptable nodes, the fixed nodes are re-added and the full list is re-sorted as detailed in \autoref{section2.1.2}.

A key challenge is labeling Dirichlet boundary nodes (see $\Lambda_{D}$ in \autoref{section2.1.2}), especially in complex cases like the L-shape domain or in the presence of interior interfaces, where node classification must be updated dynamically based on position. Standard conditional logic (namely, \texttt{if} statements) for this dynamic labeling is incompatible with JAX's JIT compiler. To address this, we employ JIT- and AD-compatible JAX-native functions from the \texttt{jax.numpy} module, like \texttt{jax.numpy.sign}, to create masks that identify nodes based on their coordinates. For example, to determine if $x\geq\alpha$ one can use $\chi_{[\alpha,\infty)}(x)=\mathrm{sgn}\big(1-\tfrac{1}{2}\big(1+\mathrm{sgn}(-(x-\alpha))\big)\big)$ as an indicator function.

\subsection{Assembly, element routine, and boundary conditions}

While JAX supports sparse matrix formats like Coordinate format (COO) and Compressed Sparse Row (CSR), custom routines were needed for some format conversions. In addition, implementing Dirichlet boundary conditions presented challenges due, in part, to the difficulty in isolating submatrices of the sparse stiffness matrix without using JIT-incompatible conditional statements. The solution involved modifying the element stiffness matrix directly (adding auxiliary rows and columns) during the element routine. Existing Dirichlet nodes were then carefully treated by saving the necessary information for inclusion in the force vector. This approach avoided conditional statements and dynamically-sized arrays to improve efficiency and leverage the advantages of JIT compilation, and produced a sparse stiffness matrix with the desired structure while facilitating the incorporation of the relevant Dirichlet and Neumann boundary condition data during the assembly of the force vector.

Material properties were assumed constant per element but could be heterogeneous across the domain. These were handled using JIT-compatible techniques via sign function tricks to produce indicator functions that facilitate spatial localization, as described above in \autoref{sec:Meshgenerationandlabeling}. Element integrals of shape functions in the stiffness matrix, which only involve polynomials and piecewise constant material properties, were computed exactly using numerical quadrature. However, element integrals in the force vector sometimes involve integrands that are \emph{not} polynomial, which, as described in \autoref{sec:quadraturevsexactintegration}, can lead to major issues. In these cases, for our examples, we implemented specialized routines that evaluate the integrals analytically, thus ensuring their exact evaluation. The one exception was the example in \autoref{ss:exp:doublesigmoid}, where we used Gauss-Legendre quadrature with $50^2$ integration points per element.

\subsection{Solvers}
\label{sec:solverspecifics}

We used sparse linear algebra routines suitable for the matrices arising in FEM to solve the linear systems $\bs{\mathsf{B}}_\theta \bs{\mathsf{c}}_\theta = \bsfell_\theta$.
We considered the \texttt{jax.experimental.sparse.linalg} library, noting its primary support for CUDA GPU backends at the time of implementation, with a fallback to \texttt{scipy.sparse.linalg.spsolve} for CPU execution.
These sparse solver routines are compatible with JAX's AD but lack support for batched operations (e.g., \texttt{vmap}). This may deteriorate performance when considering large batches of parameter samples during the training phase.
However, as noted in \autoref{remark:vanishing_partial}, it is not necessary to automatically differentiate with respect to $\bs{\mathsf{c}}_\theta$, and this can be implemented in JAX by overriding its default AD behavior via careful usage of the \texttt{@jax.custom\_vjp} decorator. Thus, if large batches are being considered or the aforementioned solvers are unsuitable, one may exploit JAX's batched CSR matrix format features and replace the solver by a more efficient (possibly iterative) one to gain efficiency. In practice, this includes the use of external high-performance solvers such as PETSc \cite{petsc1} or MUMPS \cite{MUMPS1}, which can significantly increase the computational efficiency in large-scale settings.

\section{Numerical results}\label{section4}

\subsection{Benchmark problems, training procedure, and relative errors}

\subsubsection{Benchmark problems}

We consider two elliptic boundary-value model problems, distinguished by the parameter dependence residing in either the left- or right-hand side of Poisson's equation:
\begin{subequations}\label{eq:ellipticmodelproblem}
\begin{equation}\label{eq:model_problem1}
\mathopen{}\left\lbrace
\begin{alignedat}{4}
-\Delta \uexactpar & = f^\sigma, &&\quad \text{in }\Omega, \\
\uexactpar &= 0, &&\quad \text{on }\Gamma_{D},\\
\nabla \uexactpar \cdot n & = g^\sigma, &&\quad \text{on }\Gamma_{N},
\end{alignedat}
\right.
\end{equation}
\begin{equation}\label{eq:model_problem2}
\mathopen{}\left\lbrace
\begin{alignedat}{4}
-\nabla \cdot \big(\sigma \nabla \uexactpar\big) & = f, &&\quad \text{in } \Omega, \\
\uexactpar &= 0, &&\quad \text{on }\Gamma_{D},\\
\sigma\nabla \uexactpar \cdot n & = g, &&\quad \text{on }\Gamma_{N},
\end{alignedat}
\right.
\end{equation}
where $\Omega$ is the spatial domain, $\Gamma_D$ and $\Gamma_N$ are the Dirichlet and Neumann boundaries, respectively, $g\in L^2(\Gamma_N)$ is the Neumann condition, and $n$ is the outward unit normal vector on $\Gamma_N$. Regarding the parametric behavior, $f^\sigma\in L^2(\Omega)$ is $\sigma$-dependent in \eqref{eq:model_problem1}, while $0<\sigma\in L^\infty(\Omega)$ is piecewise constant and $f\in L^2(\Omega)$ is fixed in \eqref{eq:model_problem2}.
\end{subequations}

\begin{subequations}
The variational formulations of \eqref{eq:ellipticmodelproblem} then read as: seek $u^\sigma\in H^1_D(\Omega) = \{v\in H^1(\Omega) : v{\vert_{\Gamma_D}} = 0\}$ such that
\begin{alignat}{5}
b(u^\sigma,v) &:=\int_\Omega \nabla u^\sigma\cdot\nabla v\,\mathrm{d}\Omega &&= \int_\Omega f^\sigma v\,\mathrm{d}\Omega + \int_{\Gamma_N} g^\sigma v\,\mathrm{d}\Gamma =: \ell^\sigma(v), &&\qquad \forall v\in H^1_D(\Omega).\\
b^\sigma(u^\sigma,v) &:=\int_\Omega \sigma\nabla u^\sigma\cdot\nabla v\,\mathrm{d}\Omega &&= \int_\Omega f\phantom{^\sigma}v\,\mathrm{d}\Omega + \int_{\Gamma_N} g\phantom{^\sigma}v\,\mathrm{d}\Gamma =: \phantom{^\sigma}\ell(v), &&\qquad \forall v\in H^1_D(\Omega).
\end{alignat}
\end{subequations} 
We highlight that the $\Vert\cdot\Vert_b$ norm is $\sigma$-independent in \eqref{eq:model_problem1}, i.e., $\Vert\cdot\Vert_b^2 = b(\cdot,\cdot)=\big(\nabla(\cdot),\nabla(\cdot)\big)_{L^2(\Omega)}=|\cdot|_{H^1(\Omega)}^2$, while it is $\sigma$-dependent in \eqref{eq:model_problem2}, i.e., $\Vert\cdot\Vert_{b^\sigma}^2 = b^\sigma(\cdot,\cdot)=\big(\sigma\nabla(\cdot),\nabla(\cdot)\big)_{L^2(\Omega)}=\|\sqrt{\sigma}\nabla(\cdot)\|_{L^2(\Omega)}^2$.

In what follows, $\Omega\subseteq\mathbb{R}^d$, where $d$ is the spatial dimension. 
We will refer to the exact solution as $u=u^\sigma$, whereas $u_\theta$ denotes the $r$-adaptive solution obtained via our $\theta$-based discretization for non-parametric problems (\autoref{section2.1.2}) trained at a specific value of $\sigma$, $u_\theta^\sigma$ is the $r$-adaptive solution at $\sigma$ using our NN-based discretization for parametric problems (\autoref{section2.2.2}), and $u_h = u^\sigma_h$ indicates the FEM solution on a uniform mesh of size $h=\tfrac{1}{N}$.

\subsubsection{Neural network architecture and training procedure}
\label{sec:trainingprocedure}

For the case of parametric NNs, the parameter space $\Sigma$ often comprises multiple individual parameters. Mathematically, this space is represented as a Cartesian product $\Sigma = \Sigma_1 \times \dots \times \Sigma_P$. Consequently, any specific parameter configuration $\sigma$ within this space is a tuple $\sigma = (s_1, \dots, s_P)$. When we generate a finite sample set $\{\sigma_1, \dots, \sigma_J\} \subset \Sigma$, such as the one used in \eqref{eq:loss_parametric}, each sample tuple $\smash{\sigma_i = (s_1^{(i)}, \dots, s_P^{(i)})}$ is formed by drawing its $p$-th component, $\smash{s_p^{(i)}}$, from the corresponding component space $\Sigma_p$.

For our parameter data, we sampled each $\Sigma_p$ and constructed a set of tuples $\sigma\in\Sigma$ as a `grid' of all possible combinations. 
Each space $\Sigma_p$ takes values in a compact interval, and was carefully selected, especially when the relevant parameter led to singular or localized solution behavior close to an interval extreme.
In these cases, the associated distribution was often biased toward that extreme, to improve the accuracy of the trained parametric NN in that vicinity (see the examples below for the specific distributions used).
Moreover, with the same purpose, the `corners' in the grid of parameter tuples, i.e., $\sigma=(s_1,\ldots,s_P)$, where each $s_p$ is an interval extreme, were explicitly selected as part of the training set $\{\sigma_1,\sigma_2,\ldots,\sigma_J\}\subset\Sigma$. 
The training set, which contained 70\% of the grid tuples, was otherwise randomly selected, leaving behind 30\% of grid tuples as the test set.

We used the Adam optimizer \cite{kingma2014adam} for training with learning rates (LR) in the range of $10^{-5}$ to $10^{-2}$, and trained the parametric NNs carefully until we observed convergence, which was monitored using a fixed test subset of ten randomly chosen tuples from the test set.
This sometimes required an epoch-dependent LR-strategy, where an epoch is a sweep over the entire training set.
The specific LR-strategy is detailed in each example below.

Regarding the NN architecture used for the parametric NN examples, the input layer size is $P$, which corresponds to the $\sigma$-tuple size, and it is followed by two densely connected layers with 10 nodes each, and a final output layer of the size of the number of $r$-adaptable nodes (see the size of $\delta(\sigma)$ in \eqref{eq:softmax}).
We used LeCun initialization \cite{lecun2002efficient}, and hyperbolic tangent as activation functions (i.e., $\varphi$ in \eqref{eq:net2}) at every neuron.
Lastly, for the non-parametric NN described in \autoref{section2.1.2}, initial weights were selected as $\theta^{(0)}=(0,\ldots,0)$.

\subsubsection{Relative errors}
\label{sec:relativeerrors}

For the non-parametric examples, we report the relative errors using the $\|\cdot\|_b$ norm. Note that, in view of \eqref{eq:Ritzidentity}, $\mathcal{J}(u_\theta)-\mathcal{J}(u)=\tfrac{1}{2}\|u_\theta-u\|_b^2$, and that $\ell(u)=b(u,u)=\|u\|_b^2$, so $\mathcal{J}(u)=\tfrac{1}{2}b(u,u)-\ell(u)=-\tfrac{1}{2}\|u\|_b^2$. 
Thus, we define the relative errors for the $r$-adaptive mesh and the corresponding uniform grid as
\begin{subequations}
\begin{align}
    e_\theta &:=\Big(\frac{\mathcal{J}(u)-\mathcal{J}(u_\theta)}
        {\mathcal{J}(u)}\Big)^{\tfrac{1}{2}}
        =\frac{\Vert u - u_\theta\Vert_b}{\Vert u \Vert_b},\\     
    e_h &:=\Big(\frac{\mathcal{J}(u)-\mathcal{J}(u_h)}
        {\mathcal{J}(u)}\Big)^{\tfrac{1}{2}}
        = \frac{\Vert u - u_h\Vert_b}{\Vert u \Vert_b}.
\end{align}
\end{subequations}
For the parametric examples, we proceed analogously, noting that on this occasion we can use either the $\sigma$-balanced or usual Ritz energy for the $\sigma$-specific relative error computation,
\begin{subequations}
\begin{align}
    e_\theta^\sigma 
        &:= \Big(\frac{\tilde{\mathcal{J}}^\sigma(u^\sigma)-\tilde{\mathcal{J}}^\sigma(u^\sigma_\theta)}{\tilde{\mathcal{J}}^\sigma(u^\sigma)} \Big)^{\tfrac{1}{2}}
        = \Big(\frac{\mathcal{J}^\sigma(u^\sigma)-\mathcal{J}^\sigma(u^\sigma_\theta)}{\mathcal{J}^\sigma(u^\sigma)}\Big)^{\tfrac{1}{2}} 
        = \frac{\Vert u^\sigma - u_\theta^\sigma\Vert_{b^\sigma}}{\Vert u^\sigma \Vert_{b^\sigma}},\\         
    e_h^\sigma 
        &:= \Big(\frac{\tilde{\mathcal{J}}^\sigma(u^\sigma)-\tilde{\mathcal{J}}^\sigma(u^\sigma_h)}{\tilde{\mathcal{J}}^\sigma(u^\sigma)} \Big)^{\tfrac{1}{2}}
        = \Big(\frac{\mathcal{J}^\sigma(u^\sigma)-\mathcal{J}^\sigma(u^\sigma_h)}{\mathcal{J}^\sigma(u^\sigma)}\Big)^{\tfrac{1}{2}}
        = \frac{\Vert u^\sigma - u_h^\sigma\Vert_{b^\sigma}}{\Vert u^\sigma \Vert_{b^\sigma}}.
\end{align}
\end{subequations} 
The mean and maximum of these errors are then computed over the $\sigma$'s in the training and test datasets.
Lastly, as the optimizer iterated, we used $e_{\theta}$ to monitor the convergence of the non-parametric case described in \autoref{section2.1.2}, whereas for the parametric case we used the average of $e_\theta^\sigma$ over the ten $\sigma$-tuples in a fixed test subset $\{\sigma_t\}_{t=1}^{10}$,
\begin{equation}
e_\theta^{\text{test}} := \frac{1}{10}\sum_{t=1}^{10} e_\theta^{\sigma_t}\,.
\end{equation} 

We emphasize that exact integration is critical for the above equalities to hold, as they involve Ritz energy computations, especially when the load term $\ell$ is challenging to integrate exactly (recall \autoref{sec:quadraturevsexactintegration}).

\subsection{One-dimensional experiments}\label{ss:onedimensionalexamples}

We present three one-dimensional benchmark problems, in both their non-parametric and parametric versions. In all examples, the domain is $\Omega = (0, 1)$. 
In light of the comments in \autoref{sec:quadraturevsexactintegration}, all integrals were computed exactly, in some cases via analytic integration routines specialized to the given problem.
We utilized the Adam optimizer in all of our experiments, both in the non-parametric and parametric cases.
The linear system is solved using \texttt{jax.experimental.sparse.linalg} within the JAX AD framework in the non-parametric examples, while we avoided differentiating with respect to the solution and used \texttt{scipy.sparse.linalg.spsolve} in the parametric ones, as detailed in \autoref{sec:solverspecifics}, using batches of size $10$.

\subsubsection{Arctangent sigmoid functions}
\label{ss:highgradient1d}

We begin by considering the model problem~\eqref{eq:model_problem1}, using an arctangent sigmoid function as a manufactured solution that features a sharp internal gradient given by $u(x) = \arctan(\alpha(x-s)) + \arctan(\alpha s)$. The associated forcing function is $f(x) = -u''(x)=\tfrac{2\alpha^3(x-s)}{(1+\alpha^2 (x-s)^2)^2}$. We consider a homogeneous Dirichlet boundary condition at $\Gamma_D=\{0\}$ and a Neumann boundary condition at $\Gamma_N = \{1\}$, so that $u(0)=0$ and $u'(1)=g:=\tfrac{\alpha}{1+\alpha^2(1-s)^2}$.

\paragraph{Non-parametric case.}
We fix the parameters $\bs{\sigma}=(\alpha,s):=(10, 0.5)$.
\autoref{fig:exampleI2aerror}(a) shows the evolution of the relative error during the $r$-adaptive optimization process for different numbers of elements, $N$, trained at a fixed LR of $10^{-2}$. 
In all cases, the relative error stagnates fairly quickly, at about 300 iterations, with only marginal improvements when employing further iterations. This indicates a rather quick convergence towards a (possibly locally) optimal $r$-adaptive grid. 
The errors for the $r$-adapted mesh and uniform mesh, for different values of $N$, are displayed in \autoref{fig:exampleI2aerror}(b), where the expected convergence rates are observed in both cases, but the $r$-adaptive method exhibits a better multiplicative constant than when considering the uniform mesh approach. 

\begin{figure}[!htb]
    \centering
    \begin{subfigure}[t]{0.4\textwidth}
        \resizebox{\textwidth}{!}{%
        \input{figures/example_I_2a/error-epochs-all}
        }
        \caption{Relative error during $r$-adaptive training.}
        \label{fig:exampleI2a-errors-epochs-multipleN}
    \end{subfigure}\mbox{}\hspace{4mm}
    \begin{subfigure}[t]{0.4\textwidth}
        \resizebox{\textwidth}{!}{%
        \input{figures/example_I_2a/error_history}%
        }
        \caption{Relative error convergence vs.~$N$.}
        \label{fig:exampleI2a-convergence-32elements}
    \end{subfigure}
    \caption{Training histories during $r$-adaptive optimization for various mesh sizes $N$ (left panel) and relative errors of the uniform and $r$-adapted finite element solutions as a function of $N$ (right panel) for the non-parametric model problem \eqref{eq:model_problem1} with manufactured solution $u(x) = \arctan(10(x-0.5)) + \arctan(5)$.\vspace{-2.1mm}
    }
    \label{fig:exampleI2aerror}
\end{figure}

In particular, Figures \ref{fig:exampleI2aparticular}(a) and \ref{fig:exampleI2aparticular}(b) compare FEM approximations using a uniform mesh versus an $r$-adaptive mesh using $N=32$ elements. The mesh generated by the $r$-adaptive procedure clusters the nodes around the regions with rapid solution variations, ensuring a superior approximation. This demonstrates the method's capability to adjust the mesh density according to solution features during training.

\begin{figure}[htb!]
    \centering
    \begin{subfigure}[t]{0.33\textwidth}
        \resizebox{\textwidth}{!}{%
        \input{figures/example_I_2a/FEM-32-element}%
        }
        \caption{Uniform mesh solution.
        }
        \label{fig:exampleI2a-uniform-32elements}
    \end{subfigure}\mbox{}\hspace{6mm}
    \begin{subfigure}[t]{0.33\textwidth}
        \resizebox{\textwidth}{!}{%
        \input{figures/example_I_2a/rAdaptive-32-element}%
        }
        \caption{$r$-Adapted mesh solution.
        }
        \label{fig:exampleI2a-rAdaptive-32elements}
    \end{subfigure}
    \caption{Finite element solutions using a uniform and $r$-adapted mesh according to the results in \autoref{fig:exampleI2aerror} for $N=32$.
    \vspace{-4mm}
    }
    \label{fig:exampleI2aparticular}
\end{figure}

\paragraph{Parametric case.} We fix the number of elements to $N$ and use the Ritz energy of the solution at a uniform mesh of size $h=\tfrac{1}{N}$ to compute the balanced Ritz functional in \eqref{eq:balancedRitzdef}. To bias the NN towards high-gradient solutions associated with high values of $\alpha$, we use $100$ data points coming from a reversed base-2 logarithmic distribution, ranging from $\alpha_{\min}=1$ to $\alpha_{\max}=50$, i.e., $\alpha_j=\alpha_{\min}+\alpha_{\max}-2^{\beta_j}$ where the $\beta_j$ are equispaced in $[\log_2(\alpha_{\min}),\log_2(\alpha_{\max})]$. 
For the $s$ parameter, we use $100$ points equidistant in $[0.2,0.8]$. 
Then, the $100^2$ combinations of all $(\alpha,s)$ tuples were separated into training and test sets, as described in \autoref{sec:trainingprocedure}, ensuring the parameter `corners', i.e., $(\alpha, s)\in\{(1, 0.2), (1, 0.8), (50, 0.2), (50, 0.8)\}$, were part of the training set.
Thus, there were $7\mathord{,}000$ tuples in the training set, with a batch of size $10$ being processed per optimizer iteration, so that each epoch consisted of 700 iterations. 
The LR was set to $10^{-2}$ for 20 epochs, and then to $10^{-3}$ from then onward (until 50 epochs). 
The results of the training for $N=16$ and $N=256$ are shown in \autoref{fig:exampleI2b-training-errors}.

\begin{figure}[htb!]
    \centering
    \begin{subfigure}[t]{0.4\textwidth}
        \resizebox{\textwidth}{!}{%
        \input{figures/example_I_2b/training-error-N16}
        }
        \caption{Case $N=16$.}
        \label{fig:exampleI2b-training-errors-N16}
    \end{subfigure}\mbox{}\hspace{4mm}
    \begin{subfigure}[t]{0.4\textwidth}
        \resizebox{\textwidth}{!}{%
        \input{figures/example_I_2b/training-error-N256}%
        }
        \caption{Case $N=256$.}
        \label{fig:exampleI2b-training-errors-N256}
    \end{subfigure}
    \caption{Training histories during $r$-adaptivity for $N=16$ and $N=256$ in the parametric model problem \eqref{eq:model_problem1} with manufactured solution $u^{\bs{\sigma}}(x)=\arctan(\alpha(x-s))+\arctan(\alpha s)$ with $\bs{\sigma}=(\alpha,s)$.
    }
    \label{fig:exampleI2b-training-errors}
\end{figure}

\autoref{tab:errors_sigmoid1d} compares the performance between the $r$-adaptive method and a uniform grid at different values of $N$. Focusing on the test set and $N=16$, the $r$-adaptive approach achieves a mean relative error of $0.084$ and a maximum relative error of $0.096$. In contrast, the uniform grid results in significantly higher errors, with a mean of $0.38$ and a maximum of $0.58$. The training set exhibits similar improvements, and so do the errors for $N=256$. These findings highlight the substantial accuracy improvements (approximately a factor of 5) provided by the parametric $r$-adaptive method over a standard uniform mesh across the parameter space. Furthermore, the consistent reduction of the relative errors observed on the test set confirms the NN's ability to effectively generate suitable adapted meshes, even for parameter combinations not encountered during training.

\begin{table}[htb!]
    \centering
    \caption{Relative error statistics for the uniform and the parametric $r$-adapted solutions associated with \autoref{fig:exampleI2b-training-errors}.}
    \label{tab:errors_sigmoid1d}
    \pgfplotstableread[col sep=comma]{data/example_I_2b/Error_result.csv}\dataNThirtyTwo
    \pgfplotstableread[col sep=comma]{data/example_I_2b/Error_result_N256.csv}\dataNTwoFiftySix

    \begin{tabular}{c c c c c c}
        \toprule
        $N$ & Dataset &
        \makecell{$r$-adaptive \\ mean $e_\theta^{\sigma}$} &
        \makecell{$r$-adaptive \\ max $e_\theta^{\sigma}$} &
        \makecell{uniform grid \\ mean $e_h^{\sigma}$} &
        \makecell{uniform grid \\ max $e_h^{\sigma}$} \\
        \midrule

        \multirow{2}{*}{16}
        & \pgfplotstablegetelem{0}{Set}\of\dataNThirtyTwo
        \pgfplotsretval &
        \pgfplotstablegetelem{0}{Mean Relative Error}\of\dataNThirtyTwo
        \pgfmathprintnumber[fixed, zerofill, precision=4]{\pgfplotsretval} &
        \pgfplotstablegetelem{0}{Max Relative Error}\of\dataNThirtyTwo
        \pgfmathprintnumber[fixed, zerofill, precision=4]{\pgfplotsretval} &
        \pgfplotstablegetelem{0}{Uniform grid Mean Relative error}\of\dataNThirtyTwo
        \pgfmathprintnumber[fixed, zerofill, precision=4]{\pgfplotsretval} &
        \pgfplotstablegetelem{0}{Uniform grid Max Relative error}\of\dataNThirtyTwo
        \pgfmathprintnumber[fixed, zerofill, precision=4]{\pgfplotsretval} \\
        & \pgfplotstablegetelem{1}{Set}\of\dataNThirtyTwo
        \pgfplotsretval &
        \pgfplotstablegetelem{1}{Mean Relative Error}\of\dataNThirtyTwo
        \pgfmathprintnumber[fixed, zerofill, precision=4]{\pgfplotsretval} &
        \pgfplotstablegetelem{1}{Max Relative Error}\of\dataNThirtyTwo
        \pgfmathprintnumber[fixed, zerofill, precision=4]{\pgfplotsretval} &
        \pgfplotstablegetelem{1}{Uniform grid Mean Relative error}\of\dataNThirtyTwo
        \pgfmathprintnumber[fixed, zerofill, precision=4]{\pgfplotsretval} &
        \pgfplotstablegetelem{1}{Uniform grid Max Relative error}\of\dataNThirtyTwo
        \pgfmathprintnumber[fixed, zerofill, precision=4]{\pgfplotsretval} \\
        \midrule

        \multirow{2}{*}{256}
        & \pgfplotstablegetelem{0}{Set}\of\dataNTwoFiftySix
        \pgfplotsretval &
        \pgfplotstablegetelem{0}{Mean Relative Error}\of\dataNTwoFiftySix
        \pgfmathprintnumber[fixed, zerofill, precision=4]{\pgfplotsretval} &
        \pgfplotstablegetelem{0}{Max Relative Error}\of\dataNTwoFiftySix
        \pgfmathprintnumber[fixed, zerofill, precision=4]{\pgfplotsretval} &
        \pgfplotstablegetelem{0}{Uniform grid Mean Relative error}\of\dataNTwoFiftySix
        \pgfmathprintnumber[fixed, zerofill, precision=4]{\pgfplotsretval} &
        \pgfplotstablegetelem{0}{Uniform grid Max Relative error}\of\dataNTwoFiftySix
        \pgfmathprintnumber[fixed, zerofill, precision=4]{\pgfplotsretval} 
        \\
        & \pgfplotstablegetelem{1}{Set}\of\dataNTwoFiftySix
        \pgfplotsretval &
        \pgfplotstablegetelem{1}{Mean Relative Error}\of\dataNTwoFiftySix
        \pgfmathprintnumber[fixed, zerofill, precision=4]{\pgfplotsretval} &
        \pgfplotstablegetelem{1}{Max Relative Error}\of\dataNTwoFiftySix
        \pgfmathprintnumber[fixed, zerofill, precision=4]{\pgfplotsretval} &
        \pgfplotstablegetelem{1}{Uniform grid Mean Relative error}\of\dataNTwoFiftySix
        \pgfmathprintnumber[fixed, zerofill, precision=4]{\pgfplotsretval} &
        \pgfplotstablegetelem{1}{Uniform grid Max Relative error}\of\dataNTwoFiftySix
        \pgfmathprintnumber[fixed, zerofill, precision=4]{\pgfplotsretval} \\
        \bottomrule
    \end{tabular}
\end{table}

\autoref{fig:exampleI2b-test} provides qualitative results, displaying the solutions ($u$, $u_{\theta}$ and $u_{\theta}^{\sigma}$) for three distinct test parameter pairs: $(\alpha, s) \in\{(49.27, 0.42), (19.88, 0.55), (9.96, 0.76)\}$.  
The visual results in \autoref{fig:exampleI2b-test} support that the parametric NN-generated meshes (see $u_{\theta}^{\sigma}$) and the optimized non-parametric adaptive meshes (see $u_{\theta}$) are comparable, with both effectively capturing the exact solution behavior.
This confirms the approach's ability to handle parametric dependencies and adapt the discretization efficiently.

\begin{figure}[htb!]
    \centering
    \begin{subfigure}[t]{0.31\textwidth}
        \resizebox{\textwidth}{!}{%
        \input{figures/example_I_2b/test1}%
        }
        \caption{Test sample $\boldsymbol{\sigma} = (49.57, 0.42)$:\\ $e_\theta = 0.0914$ and $e_\theta^{\sigma} = 0.0877$.}
    \end{subfigure}\mbox{}\hspace{2mm}
    \begin{subfigure}[t]{0.31\textwidth}
        \resizebox{\textwidth}{!}{%
        \input{figures/example_I_2b/test2}%
        }
        \caption{Test sample $\boldsymbol{\sigma}=(19.88,0.55)$:\\ $e_\theta = 0.0803$ and $e_\theta^{\sigma} = 0.0819$.}
    \end{subfigure}\mbox{}\hspace{2mm}
    \begin{subfigure}[t]{0.31\textwidth}
        \resizebox{\textwidth}{!}{%
        \input{figures/example_I_2b/test3}%
        }
        \caption{Test sample $\boldsymbol{\sigma}=(9.96, 0.76)$:\\ $e_\theta = 0.0715$ and $e_\theta^{\sigma} = 0.0697$.}
    \end{subfigure}
    \caption{Finite element solutions using a uniform and $r$-adapted mesh on three test samples for $N=16$. In each case, we show the exact solution $u$, the solution $u_\theta$ coming from a case-by-case non-parametric $r$-adaptivity, and the solution $u^\sigma_\theta$ coming from the parametric $r$-adaptivity produced after the training in \autoref{fig:exampleI2b-training-errors}.
    }\label{fig:exampleI2b-test}
\end{figure}

Curiously, we observed that in some cases the Ritz energy associated with $u_{\theta}$ was very similar, yet slightly higher (i.e., worse), than the one produced by $u_{\theta}^{\sigma}$, despite the fact that $u_{\theta}$ corresponded to an $r$-adapted solution especially trained at a fixed parameter tuple.
We believe this happened because the non-parametric network, which is initialized at a uniform mesh, gets stuck at a local minimum different from the global minimum as it $r$-adapts. 
Meanwhile, the parametric neural network appears to be more versatile in this respect, possibly because it is initialized differently and because it has a much wider landscape of $\bs{\sigma}=(\alpha,s)$ paths accesible to reach the global minimum as it is trained.

\subsubsection{Singular power solutions}\label{ss:singularsolutions1d}
Consider the model problem \eqref{eq:model_problem1} with the manufactured solution $u(x) = x^{\sigma}$, and the boundary conditions set as $\Gamma_D=\{0\}$ and $\Gamma_N = \{1\}$, so that $-u''(x)=f(x):=\sigma(1-\sigma)x^{\sigma-2}$ with $u(0)=0$ and $u'(1)=\sigma$.

\paragraph{Non-parametric case.}
We fix $\sigma = 0.7$, which leads to an exact solution that exhibits a singularity at $\Gamma_D$, where the derivative blows up.
As before, \autoref{fig:exampleI1a-errors-epochs}(a) shows the relative error convergence during optimization, carried out for $20\mathord{,}000$ iterations at a fixed LR of $10^{-2}$ for each of the values of $N$ displayed. 
\autoref{fig:exampleI1a-errors-epochs}(b) shows the $r$-adapted relative error for each $N$ along with the error for a uniform mesh of size $h=\tfrac{1}{N}$. 
From interpolation theory, we expect the convergence under uniform refinements to behave like $\mathcal{O}(h^{\min\{1,\sigma-0.5\}})$. Indeed, that is precisely the ratio observed in the plot, which exhibits a convergence rate of about $0.2$.
Meanwhile, the $r$-adaptive approach achieves a convergence rate of $0.98$, significantly outperforming the uniform refinements, and being close to the theoretical optimal rate of $1$ (see \cite{Strouboulis} and \cite{babuvska1979direct}). 

\begin{figure}[htb!]
    \centering
    \begin{subfigure}[t]{0.4\textwidth}
        \resizebox{\textwidth}{!}{%
        \input{figures/example_I_1a/error-epochs-all}
        }
        \caption{Relative error during $r$-adaptive training.}
        \label{fig:exampleI1a-errors-epochs-multipleN}
    \end{subfigure}\mbox{}\hspace{4mm}
    \begin{subfigure}[t]{0.41\textwidth}
        \resizebox{\textwidth}{!}{%
        \input{figures/example_I_1a/error_history}%
        }
        \caption{Relative error convergence vs.~$N$.}
        \label{fig:exampleI1a-convergence-32elements}
    \end{subfigure}
    \caption{Training histories during $r$-adaptive optimization for various mesh sizes $N$ (left panel) and relative errors of the uniform and $r$-adapted finite element solutions as a function of $N$ (right panel) for the non-parametric model problem \eqref{eq:model_problem1} with manufactured solution $u(x) = x^{0.7}$.
    }
    \label{fig:exampleI1a-errors-epochs}
\end{figure}

Moreover, Figure \ref{fig:singularpowerconv} illustrates the FEM approximations using a uniform and an $r$-adapted mesh with $N=32$ elements. Notably, the $r$-adaptive solution demonstrates node concentration toward the singularity, aligning with the optimal distribution discussed in \cite{Strouboulis}.

\begin{figure}[htb!]
    \centering
    \begin{subfigure}[t]{0.33\textwidth}
        \resizebox{\textwidth}{!}{%
        \input{figures/example_I_1a/FEM-32-element}%
        }
        \caption{Uniform mesh solution.}
        \label{fig:exampleI1a-uniform-32elements}
    \end{subfigure}\mbox{}\hspace{6mm}
    \begin{subfigure}[t]{0.33\textwidth}
        \resizebox{\textwidth}{!}{%
        \input{figures/example_I_1a/rAdaptive-32-element}%
        }
        \caption{$r$-Adapted mesh solution.}
        \label{fig:exampleI1a-rAdaptive-32elements}
    \end{subfigure}
    \caption{Finite element solutions using a uniform and $r$-adapted mesh according to the results in \autoref{fig:exampleI1a-errors-epochs} for $N=32$. \vspace{-4mm}}
    \label{fig:singularpowerconv}
\end{figure}

\paragraph{Parametric case.}
Consider $200$ values of the parameter $\sigma$ taken from a base-10 logarithmic distribution ranging from $0.51$ to $5$. 
The first, second, penultimate, and last points were purposely included in the training set, which consisted of 140 values of $\sigma$, so that $14$ iterations comprised each epoch, since batches of size 10 were processed at each iteration.
\autoref{fig:exampleI1b-training-errors} shows the training history for meshes of sizes $N = 10$ and $N=256$. For $N=10$, we fixed the LR to $10^{-2}$ for $500$ epochs (i.e., $7\mathord{,}000$  iterations). For $N=256$, we set an initial LR to $10^{-2}$, which we reduced successively, every 3 epochs (i.e., 42 iterations), by $10^{-3}$, until the LR attained a value of $3\cdot{10}^{-3}$; then, at $24$ epochs ($336$ iterations), it was further adjusted to $5\cdot10^{-4}$, and at $30$ epochs ($420$ iterations) it was decreased again to $10^{-4}$, remaining fixed thereafter (until $150$ epochs or $2\mathord{,}100$  iterations).

\begin{figure}[htb!]
    \centering
    \begin{subfigure}[t]{0.4\textwidth}
        \resizebox{\textwidth}{!}{%
        \input{figures/example_I_1b/training-error-N10}
        }
        \caption{Case $N=10$.}
        \label{fig:exampleI1b-training-errors-N16}
    \end{subfigure}\mbox{}\hspace{4mm}
    \begin{subfigure}[t]{0.4\textwidth}
        \resizebox{\textwidth}{!}{%
        \input{figures/example_I_1b/training-error-N256}%
        }
        \caption{Case $N=256$.}
        \label{fig:exampleI1b-training-errors-N256}
    \end{subfigure}
    \caption{Training histories during $r$-adaptivity for $N=10$ and $N=256$ in the parametric model problem \eqref{eq:model_problem1} with manufactured solution $u^{\sigma}(x)=x^\sigma$. \vspace{-2mm}
    }
    \label{fig:exampleI1b-training-errors}
\end{figure}

Statistics in \autoref{tab:singularsolution1d_errors} show that using a parametric $r$-adaptive method significantly improves accuracy compared to a uniform grid. At $N=10$, for the test set, the adaptive meshes have a mean relative error of $0.075$ and a maximum error of $0.58$. In contrast, uniform meshes have a mean relative error of $0.14$ and a maximum error of $0.87$. The high maximum error in both cases is due to the singular solutions near $\sigma=0.5$. Similar comparisons are seen in the training set. Much more substantial improvements are observed in the $N=256$ case, in line with the results of \autoref{fig:exampleI1a-errors-epochs}. Overall, these results confirm that the NN effectively creates well-adapted meshes for new parameters, leading to better accuracy across the parameter space.

\begin{table}[htb!]
    \centering
    \caption{Relative error statistics for the uniform and the parametric $r$-adapted solutions associated with \autoref{fig:exampleI1b-training-errors}.}
    \label{tab:singularsolution1d_errors}
    \pgfplotstableread[col sep=comma]{data/example_I_1b/Error_result.csv}\dataNThirtyTwo
    \pgfplotstableread[col sep=comma]{data/example_I_1b/Error_result_N256.csv}\dataNTwoFiftySix

    \begin{tabular}{c c c c c c}
        \toprule
        $N$ & Dataset &
        \makecell{$r$-adaptive \\ mean $e_\theta^\sigma$} &
        \makecell{$r$-adaptive \\ max $e_\theta^\sigma$} &
        \makecell{uniform grid \\ mean $e_h^\sigma$} &
        \makecell{uniform grid \\ max $e_h^\sigma$} \\
        \midrule

        \multirow{2}{*}{10}
        & \pgfplotstablegetelem{0}{Set}\of\dataNThirtyTwo
        \pgfplotsretval &
        \pgfplotstablegetelem{0}{Mean Relative Error}\of\dataNThirtyTwo
        \pgfmathprintnumber[fixed, zerofill, precision=4]{\pgfplotsretval} &
        \pgfplotstablegetelem{0}{Max Relative Error}\of\dataNThirtyTwo
        \pgfmathprintnumber[fixed, zerofill, precision=4]{\pgfplotsretval} &
        \pgfplotstablegetelem{0}{Uniform grid Mean Relative error}\of\dataNThirtyTwo
        \pgfmathprintnumber[fixed, zerofill, precision=4]{\pgfplotsretval} &
        \pgfplotstablegetelem{0}{Uniform grid Max Relative error}\of\dataNThirtyTwo
        \pgfmathprintnumber[fixed, zerofill, precision=4]{\pgfplotsretval} \\
        & \pgfplotstablegetelem{1}{Set}\of\dataNThirtyTwo
        \pgfplotsretval &
        \pgfplotstablegetelem{1}{Mean Relative Error}\of\dataNThirtyTwo
        \pgfmathprintnumber[fixed, zerofill, precision=4]{\pgfplotsretval} &
        \pgfplotstablegetelem{1}{Max Relative Error}\of\dataNThirtyTwo
        \pgfmathprintnumber[fixed, zerofill, precision=4]{\pgfplotsretval} &
        \pgfplotstablegetelem{1}{Uniform grid Mean Relative error}\of\dataNThirtyTwo
        \pgfmathprintnumber[fixed, zerofill, precision=4]{\pgfplotsretval} &
        \pgfplotstablegetelem{1}{Uniform grid Max Relative error}\of\dataNThirtyTwo
        \pgfmathprintnumber[fixed, zerofill, precision=4]{\pgfplotsretval} \\
        \midrule

        \multirow{2}{*}{256}
        & \pgfplotstablegetelem{0}{Set}\of\dataNTwoFiftySix
        \pgfplotsretval &
        \pgfplotstablegetelem{0}{Mean Relative Error}\of\dataNTwoFiftySix
        \pgfmathprintnumber[fixed, zerofill, precision=4]{\pgfplotsretval} &
        \pgfplotstablegetelem{0}{Max Relative Error}\of\dataNTwoFiftySix
        \pgfmathprintnumber[fixed, zerofill, precision=4]{\pgfplotsretval} &
        \pgfplotstablegetelem{0}{Uniform grid Mean Relative error}\of\dataNTwoFiftySix
        \pgfmathprintnumber[fixed, zerofill, precision=4]{\pgfplotsretval} &
        \pgfplotstablegetelem{0}{Uniform grid Max Relative error}\of\dataNTwoFiftySix
        \pgfmathprintnumber[fixed, zerofill, precision=4]{\pgfplotsretval} 
        \\
        & \pgfplotstablegetelem{1}{Set}\of\dataNTwoFiftySix
        \pgfplotsretval &
        \pgfplotstablegetelem{1}{Mean Relative Error}\of\dataNTwoFiftySix
        \pgfmathprintnumber[fixed, zerofill, precision=4]{\pgfplotsretval} &
        \pgfplotstablegetelem{1}{Max Relative Error}\of\dataNTwoFiftySix
        \pgfmathprintnumber[fixed, zerofill, precision=4]{\pgfplotsretval} &
        \pgfplotstablegetelem{1}{Uniform grid Mean Relative error}\of\dataNTwoFiftySix
        \pgfmathprintnumber[fixed, zerofill, precision=4]{\pgfplotsretval} &
        \pgfplotstablegetelem{1}{Uniform grid Max Relative error}\of\dataNTwoFiftySix
        \pgfmathprintnumber[fixed, zerofill, precision=4]{\pgfplotsretval} \\
        \bottomrule
    \end{tabular}
\end{table}
 
\autoref{fig:exampleI1b-test} shows the predictions of the parametric NN-based model ($u^\sigma_\theta$) versus those of the non-parametric $r$-adaptive approach ($u_\theta$) in three samples with parameters $\sigma\in\{0.52, 0.83, 1.52\}$. Both approaches show very good agreement, showcasing the efficacy of the NN-based parametric model. 
Note, however, that for $\sigma=0.52$ the $r$-adapted solutions are clearly distinct, with $u_{\theta}^{\sigma}$ yielding a better approximation than the local minimizer $u_{\theta}$, in line with the comments made in \autoref{ss:highgradient1d}.

\begin{figure}[htbp!]
    \centering
    \begin{subfigure}[t]{0.31\textwidth}
        \resizebox{\textwidth}{!}{%
        \input{figures/example_I_1b/test1}%
        }
        \caption{Test sample $\sigma=0.52$:\\ $e_\theta = 0.6402$ and $e_\theta^\sigma = 0.5781$}
    \end{subfigure}\mbox{}\hspace{2mm}
    \begin{subfigure}[t]{0.31\textwidth}
        \resizebox{\textwidth}{!}{%
        \input{figures/example_I_1b/test2}%
        }
        \caption{Test sample $\sigma=0.83$:\\ $e_\theta = 0.0367$ and $e_\theta^\sigma = 0.0366$.}
    \end{subfigure}\mbox{}\hspace{2mm}
    \begin{subfigure}[t]{0.31\textwidth}
        \resizebox{\textwidth}{!}{%
        \input{figures/example_I_1b/test3}%
        }
        \caption{Test sample $\sigma=1.52$:\\ $e_\theta = 0.0376$ and $e_\theta^\sigma = 0.0378$.}
    \end{subfigure}
    \caption{Finite element solutions using a uniform and $r$-adapted mesh on three test samples for $N=10$. In each case, we show the exact solution $u$, the solution $u_\theta$ coming from a case-by-case non-parametric $r$-adaptivity, and the solution $u^\sigma_\theta$ coming from the parametric $r$-adaptivity produced after the training in \autoref{fig:exampleI1b-training-errors}.}
    \label{fig:exampleI1b-test}
\end{figure}

\subsubsection{Two-material transmission}\label{ss:twomaterials1d}
Now, for a positive constant $\sigma>0$, we consider the other model problem \eqref{eq:model_problem2} with manufactured solution $u(x) = \sin(2\pi x)/\sigma(x)$, where the heterogeneous $\sigma(x)$ is a discontinuous function defined as $\sigma(x)=1$ if $x<0.5$ and $\sigma(x)=\sigma$ if $x\geq0.5$, so that $-u''(x)=f(x):=4\pi^2\sin(2\pi x)$.
The problem is solved with $\Gamma_D=\partial\Omega=\{0,1\}$, so $u(0)=u(1)=0$. 
Note that the solution is continuous since the discontinuity in $\sigma(x)$ occurs when $\sin(2\pi x)$ vanishes at $x=0.5$. Thus, the (weak) gradient is well defined and $u\in H_0^1(\Omega)$. In addition, notice that given positive constants $\sigma_1$ and $\sigma_2$, if $u$ is the solution with $\sigma=\tfrac{\sigma_2}{\sigma_1}$, then $\tilde{u}=\tfrac{1}{\sigma_1}u$ is the solution to the problem where the heterogeneous coefficient is instead given by $\tilde{\sigma}(x)=\sigma_1\chi_{[0,0.5)}(x)+\sigma_2\chi_{[0.5,1]}(x)$ (where $\chi_A(x)$ is the indicator function of a set $A\subset\Omega$).

\paragraph{Non-parametric case.}
We fix the value $\sigma=10$. The $r$-adaptive mesh incorporates the fixed point \(x=0.5\) to account for the material interface of the function $\sigma(x)$ (recall  \eqref{eq:nodesorting}). As a result, the $r$-adaptive scheme allows for the movement of nodes across the fixed point, adapting to regions where the function changes significantly. \autoref{fig:exampleI3a-errors-epochs} shows the relative error as the neural network is trained, for each $N$, for $10\mathord{,}000$ iterations with a LR of $10^{-2}$, as well as the error convergence for the $r$-adapted meshes and corresponding equispaced meshes. Although both display linear convergence, the $r$-adaptive method achieves a slightly smaller error magnitude, demonstrating its enhanced performance. 

\begin{figure}[htb!]
    \centering
    \begin{subfigure}[t]{0.4\textwidth}
        \resizebox{\textwidth}{!}{%
        \input{figures/example_I_3a/error-epochs-all}
        }
        \caption{Relative error during $r$-adaptive training.}
        \label{fig:exampleI3a-errors-epochs-multipleN}
    \end{subfigure}\mbox{}\hspace{4mm}
    \begin{subfigure}[t]{0.4\textwidth}
        \resizebox{\textwidth}{!}{%
        \input{figures/example_I_3a/error_history}%
        }
        \caption{Relative error convergence vs.~$N$.}
        \label{fig:exampleI3a-convergence-32elements}
    \end{subfigure}
    \caption{Training histories during $r$-adaptive optimization for various mesh sizes $N$ (left panel) and relative errors of the uniform and $r$-adapted finite element solutions as a function of $N$ (right panel) for the non-parametric model problem \eqref{eq:model_problem2} with manufactured solution $u(x) = \sin(2\pi x)/\sigma(x)$ with $\sigma(x)=1$ if $x<0.5$ and $\sigma(x)=10$ otherwise.
    }
    \label{fig:exampleI3a-errors-epochs}
\end{figure}

In addition, \autoref{fig:exampleI3a} shows the particular uniform and $r$-adapted solutions when $N=32$. Notably, we observe that the $r$-adaptive mesh asymmetrically distributes $19$ elements in $(0, 0.5)$ and $13$ elements in $(0.5, 1)$, concentrating more elements in the region where the solution presents more variation. 

\begin{figure}[htbp!]
    \centering
    \begin{subfigure}[t]{0.33\textwidth}

        \resizebox{\textwidth}{!}{%
        \input{figures/example_I_3a/FEM-32-element}%
        }
        \caption{Uniform mesh solution.}
        \label{fig:exampleI3a-uniform-32elements}
    \end{subfigure}\mbox{}\hspace{6mm}
    \begin{subfigure}[t]{0.33\textwidth}
        \resizebox{\textwidth}{!}{%
        \input{figures/example_I_3a/rAdaptive-32-element}%
        }
        \caption{$r$-Adapted mesh solution.}
        \label{fig:exampleI3a-rAdaptive-32elements}
    \end{subfigure}
    \caption{Finite element solutions using a uniform and $r$-adapted mesh according to the results in \autoref{fig:exampleI3a-errors-epochs} for $N=32$.
    \vspace{-4mm}
    }
    \label{fig:exampleI3a}
\end{figure}

\paragraph{Parametric case.}
We use a base-10 logarithmic distribution for $\sigma$ ranging from $10^{-4}$ to $10^{4}$, consisting of $1\mathord{,}000$ data points, i.e., $\sigma$ is taken from $10^{\beta_j}$, where the $\beta_j$ are equispaced in $[-4,4]$. The first, second, penultimate, and last points were explicitly included in the training set, which contained $700$ values of $\sigma$, meaning each epoch consisted of $70$ iterations (with a batch of size $10$ processed per iteration). As in the non-parametric case, the node at $x=0.5$ is fixed, so for a total of $N$ elements we only have $N-2$ movable nodes.
\autoref{fig:exampleI3b-training-errors} shows the training history for $N=12$ and $N=256$ elements, in both cases using a LR of $10^{-2}$ during the very first $30$ epochs ($2\mathord{,}100$ iterations) and of $10^{-3}$ in the last $120$ epochs ($8\mathord{,}400$ iterations).

\begin{figure}[htb!]
    \centering
    \begin{subfigure}[t]{0.4\textwidth}
        \resizebox{\textwidth}{!}{%
        \input{figures/example_I_3b/training-error-N12}
        }
        \caption{Case $N=12$.}
        \label{fig:exampleI3b-training-errors-N12}
    \end{subfigure}\mbox{}\hspace{4mm}
    \begin{subfigure}[t]{0.4\textwidth}
        \resizebox{\textwidth}{!}{%
        \input{figures/example_I_3b/training-error-N256}%
        }
        \caption{Case $N=256$.}
        \label{fig:exampleI3b-training-errors-N256}
    \end{subfigure}
    \caption{Training histories during $r$-adaptivity for $N=12$ and $N=256$ in the parametric model problem \eqref{eq:model_problem2} with manufactured solution $u^{\sigma}(x)=\sin(2\pi x)/\sigma(x)$, where $\sigma(x)=1$ if $x<0.5$ and $\sigma(x)=\sigma$ if $x\geq 0.5$.
    }
    \label{fig:exampleI3b-training-errors}
\end{figure}

The results in \autoref{tab:errors_twomaterials1d} indicate that using a parametric $r$-adaptive method leads to slightly better accuracy than a uniform grid. 
Although both methods produce small errors, at $N=12$, over the test set, the adaptive meshes achieve a mean relative error of $0.11$, while the uniform meshes have a mean relative error of $0.15$. 
Similar improvements were noted in the training set, and analogous results are observed for $N=256$, indicating that the NN effectively generates well-adapted meshes, enhancing accuracy across different parameters.

\begin{table}[htb!]
    \centering
    \caption{Relative error statistics for the uniform and the parametric $r$-adapted solutions associated with \autoref{fig:exampleI3b-training-errors}.}
    \label{tab:errors_twomaterials1d}
    \pgfplotstableread[col sep=comma]{data/example_I_3b/Error_result.csv}\dataNThirtyTwo
    \pgfplotstableread[col sep=comma]{data/example_I_3b/Error_result_N256.csv}\dataNTwoFiftySix

    \begin{tabular}{c c c c c c}
        \toprule
        $N$ & Dataset &
        \makecell{$r$-adaptive \\ mean $e_\theta^\sigma$} &
        \makecell{$r$-adaptive \\ max $e_\theta^\sigma$} &
        \makecell{uniform grid \\ mean $e_h^\sigma$} &
        \makecell{uniform grid \\ max $e_h^\sigma$} \\
        \midrule

        \multirow{2}{*}{12}
        & \pgfplotstablegetelem{0}{Set}\of\dataNThirtyTwo
        \pgfplotsretval &
        \pgfplotstablegetelem{0}{Mean Relative Error}\of\dataNThirtyTwo
        \pgfmathprintnumber[fixed, zerofill, precision=4]{\pgfplotsretval} &
        \pgfplotstablegetelem{0}{Max Relative Error}\of\dataNThirtyTwo
        \pgfmathprintnumber[fixed, zerofill, precision=4]{\pgfplotsretval} &
        \pgfplotstablegetelem{0}{Uniform grid Mean Relative error}\of\dataNThirtyTwo
        \pgfmathprintnumber[fixed, zerofill, precision=4]{\pgfplotsretval} &
        \pgfplotstablegetelem{0}{Uniform grid Max Relative error}\of\dataNThirtyTwo
        \pgfmathprintnumber[fixed, zerofill, precision=4]{\pgfplotsretval} \\
        & \pgfplotstablegetelem{1}{Set}\of\dataNThirtyTwo
        \pgfplotsretval &
        \pgfplotstablegetelem{1}{Mean Relative Error}\of\dataNThirtyTwo
        \pgfmathprintnumber[fixed, zerofill, precision=4]{\pgfplotsretval} &
        \pgfplotstablegetelem{1}{Max Relative Error}\of\dataNThirtyTwo
        \pgfmathprintnumber[fixed, zerofill, precision=4]{\pgfplotsretval} &
        \pgfplotstablegetelem{1}{Uniform grid Mean Relative error}\of\dataNThirtyTwo
        \pgfmathprintnumber[fixed, zerofill, precision=4]{\pgfplotsretval} &
        \pgfplotstablegetelem{1}{Uniform grid Max Relative error}\of\dataNThirtyTwo
        \pgfmathprintnumber[fixed, zerofill, precision=4]{\pgfplotsretval} \\
        \midrule

        \multirow{2}{*}{256}
        & \pgfplotstablegetelem{0}{Set}\of\dataNTwoFiftySix
        \pgfplotsretval &
        \pgfplotstablegetelem{0}{Mean Relative Error}\of\dataNTwoFiftySix
        \pgfmathprintnumber[fixed, zerofill, precision=4]{\pgfplotsretval} &
        \pgfplotstablegetelem{0}{Max Relative Error}\of\dataNTwoFiftySix
        \pgfmathprintnumber[fixed, zerofill, precision=4]{\pgfplotsretval} &
        \pgfplotstablegetelem{0}{Uniform grid Mean Relative error}\of\dataNTwoFiftySix
        \pgfmathprintnumber[fixed, zerofill, precision=4]{\pgfplotsretval} &
        \pgfplotstablegetelem{0}{Uniform grid Max Relative error}\of\dataNTwoFiftySix
        \pgfmathprintnumber[fixed, zerofill, precision=4]{\pgfplotsretval} 
        \\
        & \pgfplotstablegetelem{1}{Set}\of\dataNTwoFiftySix
        \pgfplotsretval &
        \pgfplotstablegetelem{1}{Mean Relative Error}\of\dataNTwoFiftySix
        \pgfmathprintnumber[fixed, zerofill, precision=4]{\pgfplotsretval} &
        \pgfplotstablegetelem{1}{Max Relative Error}\of\dataNTwoFiftySix
        \pgfmathprintnumber[fixed, zerofill, precision=4]{\pgfplotsretval} &
        \pgfplotstablegetelem{1}{Uniform grid Mean Relative error}\of\dataNTwoFiftySix
        \pgfmathprintnumber[fixed, zerofill, precision=4]{\pgfplotsretval} &
        \pgfplotstablegetelem{1}{Uniform grid Max Relative error}\of\dataNTwoFiftySix
        \pgfmathprintnumber[fixed, zerofill, precision=4]{\pgfplotsretval} \\
        \bottomrule
    \end{tabular}
\end{table}

\autoref{fig:exampleI3b-test} shows the predictions of the parametric NN-based model ($u^\sigma_\theta$) versus those of the non-parametric $r$-adaptive approach ($u_\theta$) in three samples with parameters $\sigma \in\{4.6 \cdot 10^{-4}, 8.9 \cdot 10^{2}, 5.5 \cdot 10^{-1}\}$.
Even though both approaches lead to predictions that capture the exact solution (and are symmetric in each half-domain), note that the $r$-adaptive meshes are different, since at each half they contain a different number of nodes.
As pointed out before, it turns out the solutions coming from the parametric NN actually yield better results (i.e., lower errors) than those from the non-parametric NN trained at the fixed $\sigma$.
The reasons are likely the ones stated previously.
Having said that, we further experimented with the non-parametric NN, modifying its optimization procedure, so instead of using Adam, we utilized stochastic gradient descent with a learning rate of $0.01$, momentum of $0.95$, and Nesterov acceleration.
By doing so, we were able to avoid the sub-optimal local minima and match the results of the parametric NN for the first two cases in \autoref{fig:exampleI3b-test}, but retained the same (locally minimizing) solution in the last case.

\begin{figure}[htb!]
    \centering
    \begin{subfigure}[t]{0.31\textwidth}
        \resizebox{\textwidth}{!}{%
        \input{figures/example_I_3b/test1}%
        }
        \caption{Test sample $\sigma = 4.6\cdot10^{-4}$:\\
        $e_\theta = 0.1140$ and $e_\theta^\sigma = 0.0994$.}
    \end{subfigure}\mbox{}\hspace{2mm}
    \begin{subfigure}[t]{0.31\textwidth}
        \resizebox{\textwidth}{!}{%
        \input{figures/example_I_3b/test2}%
        }
        \caption{Test sample $\sigma=8.9\cdot 10^{2}$:\\
        $e_\theta = 0.1140$ and $e_\theta^\sigma = 0.0995$.}
    \end{subfigure}\mbox{}\hspace{2mm}
    \begin{subfigure}[t]{0.31\textwidth}
        \resizebox{\textwidth}{!}{%
        \input{figures/example_I_3b/test3}%
        }
        \caption{Test sample $\sigma=5.5\cdot10^{-1}$.\\
        $e_\theta = 0.1337$ and $e_\theta^\sigma = 0.1334$.}
    \end{subfigure}
    \caption{Finite element solutions using a uniform and $r$-adapted mesh on three test samples for $N=12$. In each case, we show the exact solution $u$, the solution $u_\theta$ coming from a case-by-case non-parametric $r$-adaptivity, and the solution $u^\sigma_\theta$ coming from the parametric $r$-adaptivity produced after the training in \autoref{fig:exampleI3b-training-errors}.
    \label{fig:exampleI3b-test}}
\end{figure}

\subsection{Two-dimensional experiments}
\label{ss:twodimensialexperiments}

We present two two-dimensional benchmark problems, in both their non-parametric and parametric versions. Quadrilateral tensor-product meshes are used to discretize the domains. In the first example, we consider the domain $\Omega = (0, 1)^{2}$, whilst in the second, $\Omega$ is an L-shape domain.
The Adam optimizer is used in all of our experiments, both in the non-parametric and parametric cases, with the LR varying on a case-by-case basis. 
Lastly, in the non-parametric examples the linear systems are solved using \texttt{jax.experimental.sparse.linalg} within the JAX AD framework, whereas in the parametric examples they are solved, with batches of size $1$, via \texttt{scipy.sparse.linalg.spsolve}, avoiding differentiation with respect to the solution.

\subsubsection{Arctangent sigmoid functions}\label{ss:exp:doublesigmoid} 

The two-dimensional extension of the one-dimensional arctangent sigmoid solution from Subsection \ref{ss:highgradient1d} is obtained by setting $\Omega=(0,1)^2$ and considering the Poisson problem \eqref{eq:model_problem1} with manufactured solution $u(x,y)=u_1(x)u_2(y)$, where $u_j(t) = \arctan(\alpha(t - s_j)) + \arctan(\alpha s_{j})$ for $j=1,2$. Homogeneous Dirichlet boundary conditions are established at the left and bottom boundaries (i.e. $x=0$ or $y=0$), while Neumann boundary conditions are set at the right and top boundaries (i.e. $x=1$ or $y=1$).
The function $f=-\Delta u$ can thus be analytically calculated to be the expression
\begin{equation}
   f(x,y) = u_1''(x)u_2(y) + u_1(x)u_2''(y) \qquad\text{where}\qquad
    u_j''(t)=\frac{2\alpha^3(t-s_{j})}{(1+(\alpha(t - s_{j}))^2)^2}\,,\quad j=1,2\,.
\end{equation}
The integrals involving $f$ are approximated numerically using $50^2$ Gauss-Legendre quadrature points per element, since the resulting expressions are challenging to integrate exactly. The high number of quadrature points is important to avoid the issues described in \autoref{sec:quadraturevsexactintegration}.

\paragraph{Non-parametric case.} 
Let $\alpha = 10$ and $s_1=s_2=0.05$.
For each value of $N$, corresponding to the number of elements per side, we optimize the $r$-adaptive mesh during $2\mathord{,}500$ iterations using the Adam optimizer with a fixed LR of $10^{-2}$, with the results being those of \autoref{fig:2dsigmoid-errorepoch}(a).
\autoref{fig:2dsigmoid-errorepoch}(b), like \autoref{fig:exampleI2b-test}(b), shows the $r$-adapted and uniform-mesh (with $h=\tfrac{1}{N}$) solutions converge linearly, with the $r$-adapted solutions exhibiting a better constant.

\begin{figure}[htb!]
    \centering
    \begin{subfigure}[t]{0.4\textwidth}
        \resizebox{\textwidth}{!}{%
        \input{figures/example_II_1a/error-epoch-all}
        }
        \caption{Relative error during $r$-adaptive training.}
        \label{fig:exampleII1a-errors-epochs-multipleN}
    \end{subfigure}\mbox{}\hspace{4mm}
    \begin{subfigure}[t]{0.4\textwidth}
        \resizebox{\textwidth}{!}{%
        \input{figures/example_II_1a/error_history}%
        }
        \caption{Relative error convergence vs.~$N$.}
        \label{fig:exampleII1a-convergence}
    \end{subfigure}
    \caption{Training histories during $r$-adaptive optimization for each value of $N$ (left panel), and relative errors of the uniform and $r$-adapted finite element solutions as a function of $N$ (right panel) for the non-parametric model problem \eqref{eq:model_problem1} with manufactured solution $u(x,y)=u_1(x)u_2(y)$ where $u_j(t)=\arctan(10t-0.05)+\arctan(0.05)$ for $j=1,2$.
    }
    \label{fig:2dsigmoid-errorepoch}
\end{figure}

\autoref{fig:2dsigmoidconvresults} shows the uniform and $r$-adapted predictions featuring $N=32$ elements per side for a total of $32^2$ elements. As expected, the $r$-adaptive method focuses nodes in areas with sharp gradients, thereby enhancing the accuracy of the approximation. 

\begin{figure}[htb!]
    \centering
    \begin{subfigure}[t]{0.31\textwidth}
        \resizebox{\textwidth}{!}{%
        \input{figures/example_II_1a/FEM-32-element}%
        }
        \caption{Uniform mesh solution.}
        \label{fig:exampleII1a-uniform-32elements}
    \end{subfigure}\mbox{}\hspace{4mm}
    \begin{subfigure}[t]{0.31\textwidth}
        \resizebox{\textwidth}{!}{%
        \input{figures/example_II_1a/rAdaptive-32-element}%
        }
        \caption{$r$-Adapted mesh solution.}
        \label{fig:exampleII1a-rAdaptive-32elements}
    \end{subfigure}
    \caption{
    Finite element solutions using a uniform and $r$-adapted mesh according to the results in \autoref{fig:2dsigmoid-errorepoch} for $N=32$.
    \vspace{-4mm}}
    \label{fig:2dsigmoidconvresults}
\end{figure}

\paragraph{Parametric case.} 
Let the parameter be $\bs{\sigma}=(\alpha,s_1,s_2)$. 
We sample $\alpha$ from the same inverted base-$2$ logarithmic distribution described in \autoref{ss:highgradient1d}, ranging from $20$ to $1$ with $20$ data points, whereas $s_1$ and $s_2$ are uniformly distributed in $[0.1,0.9]$ with $10$ data points each. Thus, the training set consisted of $1\mathord{,}400$ $\bs{\sigma}$-tuples, with each epoch comprised of $1\mathord{,}400$ iterations (using a batch size of $1$). \autoref{fig:exampleII1b-training-errors} shows the training evolution with $N=9$ and $N=64$ elements per side. For $N=9$, training is carried out for $1\mathord{,}000$ epochs ($1\mathord{,}400\mathord{,}000$ iterations) at an initial LR of $10^{-2}$ for the first $150$ epochs ($210\mathord{,}000$ iterations) and $10^{-3}$ thereafter, while for $N=64$ it is trained for $500$ epochs ($700\mathord{,}000$ iterations) with a LR of $10^{-2}$ for the first $20$ epochs ($28\mathord{,}000$ iterations) and $10^{-3}$ afterward. 

\autoref{tab:errors_doublesigmoid} compares the relative errors of the approximation using  NN-based $r$-adaptive meshes and uniform meshes. 

\begin{table}[hbt!]
    \centering
    \caption{Relative error statistics for the uniform and the parametric $r$-adapted solutions associated with \autoref{fig:exampleII1b-training-errors}.
    }
    \label{tab:errors_doublesigmoid}

    \pgfplotstableread[col sep=comma]{data/example_II_1b/Error_result_N9.csv}\dataNThirtyTwo
    \pgfplotstableread[col sep=comma]{data/example_II_1b/Error_result_N64.csv}\dataNTwoFiftySix

    \begin{tabular}{c c c c c c}
        \toprule
        $N$ & Dataset &
        \makecell{$r$-adaptive \\ mean $e_\theta^\sigma$} &
        \makecell{$r$-adaptive \\ max $e_\theta^\sigma$} &
        \makecell{uniform grid \\ mean $e_h^\sigma$} &
        \makecell{uniform grid \\ max $e_h^\sigma$} \\
        \midrule

        \multirow{2}{*}{9}
        & \pgfplotstablegetelem{0}{Set}\of\dataNThirtyTwo
        \pgfplotsretval &
        \pgfplotstablegetelem{0}{Mean Relative Error}\of\dataNThirtyTwo
        \pgfmathprintnumber[fixed, zerofill, precision=4]{\pgfplotsretval} &
        \pgfplotstablegetelem{0}{Max Relative Error}\of\dataNThirtyTwo
        \pgfmathprintnumber[fixed, zerofill, precision=4]{\pgfplotsretval} &
        \pgfplotstablegetelem{0}{Uniform grid Mean Relative error}\of\dataNThirtyTwo
        \pgfmathprintnumber[fixed, zerofill, precision=4]{\pgfplotsretval} &
        \pgfplotstablegetelem{0}{Uniform grid Max Relative error}\of\dataNThirtyTwo
        \pgfmathprintnumber[fixed, zerofill, precision=4]{\pgfplotsretval} \\
        & \pgfplotstablegetelem{1}{Set}\of\dataNThirtyTwo
        \pgfplotsretval &
        \pgfplotstablegetelem{1}{Mean Relative Error}\of\dataNThirtyTwo
        \pgfmathprintnumber[fixed, zerofill, precision=4]{\pgfplotsretval} &
        \pgfplotstablegetelem{1}{Max Relative Error}\of\dataNThirtyTwo
        \pgfmathprintnumber[fixed, zerofill, precision=4]{\pgfplotsretval} &
        \pgfplotstablegetelem{1}{Uniform grid Mean Relative error}\of\dataNThirtyTwo
        \pgfmathprintnumber[fixed, zerofill, precision=4]{\pgfplotsretval} &
        \pgfplotstablegetelem{1}{Uniform grid Max Relative error}\of\dataNThirtyTwo
        \pgfmathprintnumber[fixed, zerofill, precision=4]{\pgfplotsretval} \\
        \midrule

        \multirow{2}{*}{64}
        & \pgfplotstablegetelem{0}{Set}\of\dataNTwoFiftySix
        \pgfplotsretval &
        \pgfplotstablegetelem{0}{Mean Relative Error}\of\dataNTwoFiftySix
        \pgfmathprintnumber[fixed, zerofill, precision=4]{\pgfplotsretval} &
        \pgfplotstablegetelem{0}{Max Relative Error}\of\dataNTwoFiftySix
        \pgfmathprintnumber[fixed, zerofill, precision=4]{\pgfplotsretval} &
        \pgfplotstablegetelem{0}{Uniform grid Mean Relative error}\of\dataNTwoFiftySix
        \pgfmathprintnumber[fixed, zerofill, precision=4]{\pgfplotsretval} &
        \pgfplotstablegetelem{0}{Uniform grid Max Relative error}\of\dataNTwoFiftySix
        \pgfmathprintnumber[fixed, zerofill, precision=4]{\pgfplotsretval} 
        \\
        & \pgfplotstablegetelem{1}{Set}\of\dataNTwoFiftySix
        \pgfplotsretval &
        \pgfplotstablegetelem{1}{Mean Relative Error}\of\dataNTwoFiftySix
        \pgfmathprintnumber[fixed, zerofill, precision=4]{\pgfplotsretval} &
        \pgfplotstablegetelem{1}{Max Relative Error}\of\dataNTwoFiftySix
        \pgfmathprintnumber[fixed, zerofill, precision=4]{\pgfplotsretval} &
        \pgfplotstablegetelem{1}{Uniform grid Mean Relative error}\of\dataNTwoFiftySix
        \pgfmathprintnumber[fixed, zerofill, precision=4]{\pgfplotsretval} &
        \pgfplotstablegetelem{1}{Uniform grid Max Relative error}\of\dataNTwoFiftySix
        \pgfmathprintnumber[fixed, zerofill, precision=4]{\pgfplotsretval} \\
        \bottomrule
    \end{tabular}
\end{table}
For $N=64$ and the test set, the $r$-adaptive method gives a mean relative error of 0.017 and a maximum relative error of 0.022. This outperforms the results on the uniform grid, where the mean and maximum relative errors are 0.046 and 0.064, respectively. Similar numbers are obtained for the training set. The results confirm the findings in \autoref{tab:errors_sigmoid1d}, namely, the NN-based meshes result in an improvement compared to uniform meshes.

\begin{figure}[tb!]
    \centering
    \begin{subfigure}[t]{0.4\textwidth}
        \resizebox{\textwidth}{!}{%
        \input{figures/example_II_1b/training-error-N9}
        }
        \caption{Case $N=9$.}
        \label{fig:exampleII1b-training-errors-N9}
    \end{subfigure}\mbox{}\hspace{4mm}
    \begin{subfigure}[t]{0.4\textwidth}
        \resizebox{\textwidth}{!}{%
        \input{figures/example_II_1b/training-error-N64}
        }
        \caption{Case $N=64$.}
        \label{fig:exampleII1b-training-errors-N256}
    \end{subfigure}
    \caption{
    Training histories during $r$-adaptivity for $N=9$ and $N=64$ elements per side in the parametric model problem \eqref{eq:model_problem1}, with manufactured solution $u(x,y)=u_1(x)u_2(y)$, where $u_j(t)=\arctan(\alpha(t-s_j))+\arctan(\alpha s_{j})$ for $j=1,2$.
    \vspace{-2mm}}
    \label{fig:exampleII1b-training-errors}
\end{figure}

\autoref{fig:exampleII1a-test} displays the solutions $u_\theta$ and $u_\theta^\sigma$ for selected parameter values, showing relatively good agreement between them and their respective two-dimensional non-parametric and parametric meshes.
\begin{figure}[hbt!]
    \centering
    \begin{subfigure}[t]{0.05\textwidth}
        \begin{tikzpicture}[scale=0.9]
        \node[anchor=west] at (0cm,0cm) {};
            \node[anchor=west] at (-0.5cm,7.0cm) {$u_{\theta}$};
            \node[anchor=west] at (-0.5cm,2.0cm) {$u_{\theta}^{\sigma}$};
        \end{tikzpicture}
    \end{subfigure}\mbox{}\hspace{1mm}
    \begin{subfigure}[t]{0.257\textwidth}
        \resizebox{\textwidth}{!}{
        \input{figures/example_II_1b/test1_r_NN}
        }
        \caption{Test sample $\bs{\sigma} = (17.98, 0.19, 0.28)$:\\ $e_\theta=0.1878$ and $e_\theta^{\sigma}=0.1353$.}
    \end{subfigure}\mbox{}\hspace{2mm}
    \begin{subfigure}[t]{0.25\textwidth}
        \resizebox{\textwidth}{!}{
        \input{figures/example_II_1b/test2_r_NN}
        }
        \caption{Test sample $\bs{\sigma} = (16.16, 0.54, 0.46)$:\\
        $e_\theta=0.1336$ and $e_\theta^{\sigma}=0.1345$.}
    \end{subfigure}\mbox{}\hspace{2mm}
    \begin{subfigure}[t]{0.25\textwidth}
    \resizebox{\textwidth}{!}{
        \input{figures/example_II_1b/test3_r_NN}
        }
        \caption{Test sample $\bs{\sigma} = (3.92, 0.28, 0.9)$:\\
        $e_\theta=0.0567$ and $e_\theta^{\sigma}=0.0681$.}
    \end{subfigure}
    \caption{
    Finite element solutions using a uniform and $r$-adapted mesh on three test samples for $N=9$. In each case, we show the solution $u_\theta$ coming from a case-by-case non-parametric $r$-adaptivity, and the solution $u^\sigma_\theta$ coming from the parametric $r$-adaptivity produced after the training in \autoref{fig:exampleII1b-training-errors}.
    \vspace{-5mm}}
    \label{fig:exampleII1a-test}
\end{figure}

\subsubsection{Multi-material L-shape domain} 
\label{ss:exp:Lshape}

We consider the benchmark problem \eqref{eq:model_problem1} on an L-shape domain, $\Omega=(0,1)^2\setminus\big([0.5,1]\times[0,0.5]\big)$, with $f=1$ and homogeneous Dirichlet boundary conditions, i.e. $u|_{\partial\Omega}=0$. Given positive constants $\sigma_1$ and $\sigma_2$, the heterogeneous $\sigma(\bs{x})$ is given by
\begin{equation}
    \sigma(\bs{x}) 
        = \chi_{(0,0.5)\times(0.5,1)}(\bs{x})
            +\sigma_1\chi_{(0,0.5)^2}(\bs{x})
                +\sigma_2\chi_{(0.5,1)^2}(\bs{x}) =
    \begin{cases}
    \sigma_0=1\quad&\bs{x}\in(0,0.5)\times(0.5,1),\\
    \sigma_1 \quad&\bs{x}\in(0,0.5)^2,\\
    \sigma_2 \quad&\bs{x}\in(0.5,1)^2,
    \end{cases}
    \label{eq:materialproLshape}
\end{equation}
where $\chi_A(\bs{x})$ is the indicator function for the set $A\subset\Omega$.
Note that given a triplet of positive numbers $(\tilde{\sigma}_0,\tilde{\sigma}_1,\tilde{\sigma}_2)$, if $u$ is the solution of the problem above with $\sigma_1=\tfrac{\tilde{\sigma}_1}{\tilde{\sigma}_0}$ and $\sigma_2=\tfrac{\tilde{\sigma}_2}{\tilde{\sigma}_0}$, then $\tilde{u}=\tfrac{1}{\tilde{\sigma_0}}u$ is the solution of the model problem for a general three-material $\tilde{\sigma}(\bs{x})$ given by 
\begin{equation*}
    \tilde{\sigma}(\bs{x})=\tilde{\sigma}_0\sigma(\bs{x})
        =\tilde{\sigma}_0\chi_{(0,0.5)\times(0.5,1)}(\bs{x})
            +\tilde{\sigma}_1\chi_{(0,0.5)^2}(\bs{x})
                +\tilde{\sigma}_2\chi_{(0.5,1)^2}(\bs{x})\,.
\end{equation*}
Thus, with only two parameters, i.e., $\bs{\sigma}=(\sigma_1,\sigma_2)$, one is able to characterize the three-material case. 

With this $f$ and boundary conditions, the solutions are not known analytically, even for the simplest $\bs{\sigma}=(\sigma_1,\sigma_2)$.
Nevertheless, they are known to have a singular gradient at the re-entrant corner, and the `intensity' of this singularity varies with $\bs{\sigma}$. Indeed, in the limit $\sigma_2\to\infty$, leaving $\sigma_1$ constant, the singularity will disappear entirely, effectively leading to the equation being satisfied in a rectangle.
As a replacement for the exact solutions, we resort to highly accurate approximations of their Ritz energy computed in carefully refined triangular meshes with the open source FEM software \mbox{Netgen/NGSolve}\cite{schoberl1997netgen,schoberl2014c++}.
Note that this `exact' Ritz energy is enough to compute the relevant $\|\cdot\|_b$ relative errors in \autoref{sec:relativeerrors} (rather than interpolate between meshes).

Regarding our implementation, we use quadrilateral meshes with $N$ elements per side, where, to ensure the presence of the re-entrant corner in the mesh, nodes with $x=0.5$ and $y=0.5$ coordinates are fixed.
To emulate the L-shape domain, the nodes in the bottom right quadrant are tagged as zero Dirichlet nodes, and these are dynamically labeled after each $r$-adaptive step, as described in \autoref{sec:Meshgenerationandlabeling}.
Due to the simplicity of the forcing and boundary conditions (i.e., $f=1$ and $u|_{\partial\Omega}=0$), all integrals in the finite element variational formulation are computed exactly using numerical quadrature (with $2^2$ Gauss-Legendre quadrature points).

\paragraph{Non-parametric case.}
Let $\bs{\sigma}=(1,1)$, so that we solve the Poisson equation $\Delta u=1$ with $u|_{\partial\Omega}=0$. This is \emph{not} the usual L-shape domain example where the solution (i.e., $u(r,\theta)=r^{2/3}\sin(\tfrac{2}{3}\theta)$ solving $\Delta u=0$) is in $H^{1+s-\delta}(\Omega)$ for $s=\tfrac{2}{3}$ and $\delta>0$. However, it is still a classical example of a `nice' second-order operator with smooth right-hand side and boundary conditions that does not lead to a smooth solution.
Indeed, due to the nonconvex re-entrant corner, the solution has a singular gradient and does not belong to $C^1(\Omega)$. Instead, it is in $H^{1+s}(\Omega)$ where $s$ appears to be around $0.72$ as computed with Netgen/NGSolve under high-order uniform refinements.
The Ritz energy used as a reference for the relative error computations is $\mathcal{J}(u)=-0.00668986$, which was computed with order-five polynomials on an unstructured triangular mesh specially refined at the re-entrant corner using Netgen/NGSolve.

For each given (even) number of elements per side, $N$, note that the number of `active' elements in a uniform grid of size $h=\tfrac{1}{N}$ is actually $\tfrac{3}{4}N^2$ and the number of degrees of freedom is $\tfrac{3}{4}N^2-2N+1$.
At different $N$, we optimize the $r$-adaptive mesh during $100\mathord{,}000$ iterations using the Adam optimizer with a fixed LR of $10^{-2}$, as shown in \autoref{fig:2dLshape-errorepoch}(a).
The convergence in $N$ of the uniform mesh and $r$-adapted solutions is shown in \autoref{fig:2dLshape-errorepoch}(b), where the corresponding number of degrees of freedom is also $\tfrac{3}{4}N^2-2N+1$, because, in all cases observed, the nodes never crossed the midpoint.
We observe that the convergence rate with respect to $N=\tfrac{1}{h}$ under uniform refinements at $N=2^8$ is around 0.78 and appears to be decreasing (indeed, it will eventually decrease to around 0.72). Remarkably, in the $r$-adaptive case, the mean convergence rate is $1.00$, which is the best possible asymptotic convergence rate (see \cite{babuvska1979direct}).

\begin{figure}[htb!]
    \centering
    \begin{subfigure}[t]{0.4\textwidth}
        \resizebox{\textwidth}{!}{%
        \input{figures/example_II_2a/error-epoch-all}
        }
        \caption{Relative error during $r$-adaptive training.}
        \label{fig:exampleII2a-errors-epochs-multipleN}
    \end{subfigure}\mbox{}\hspace{4mm}
    \begin{subfigure}[t]{0.4\textwidth}
        \resizebox{0.98\textwidth}{!}{%
        \input{figures/example_II_2a/error_history}%
        }
        \caption{Relative error convergence vs.~$N$.}
        \label{fig:exampleII2a-convergence-32elements}
    \end{subfigure}
    \caption{
    Training histories during $r$-adaptive optimization for each value of $N$ (left panel), and relative errors of the uniform and $r$-adapted finite element solutions as a function of $N$ (right panel) for the non-parametric model problem $\Delta u=1$ in the L-shape domain with $u|_{\partial\Omega}=0$.
    }
    \label{fig:2dLshape-errorepoch}
\end{figure}

\autoref{fig:LshapeNonParametric} displays the uniform and $r$-adapted solutions for $N=32$ elements per side. The $r$-adaptive mesh is refining strongly toward the re-entrant corner and moderately toward the exterior boundaries, and the scale of the gradient magnitude clearly shows that it captures the singularity more accurately than the uniform mesh.

\begin{figure}[htb!]
    \centering
    \begin{subfigure}[t]{0.31\textwidth}
        \resizebox{\textwidth}{!}{%
        \input{figures/example_II_2a/new-FEM-32-element}%
        }
        \caption{Uniform mesh solution (top) and its gradient magnitude (bottom).}
        \label{fig:exampleII2a-uniform-32elements}
    \end{subfigure}\mbox{}\hspace{4mm}
    \begin{subfigure}[t]{0.31\textwidth}
        \resizebox{\textwidth}{!}{%
        \input{figures/example_II_2a/new-rAdaptive-32-element}%
        }
        \caption{$r$-Adapted mesh solution (top) and its gradient magnitude (bottom).}
        \label{fig:exampleII2a-rAdaptive-32elements}
    \end{subfigure}
    \caption{
    Finite element solutions using a uniform and $r$-adapted mesh according to the results in \autoref{fig:2dLshape-errorepoch} for $N=32$.
    \vspace{-4mm}}
    \label{fig:LshapeNonParametric}
\end{figure}

\paragraph{Parametric case.}
Let the parameter be $\bs{\sigma}=(\sigma_1,\sigma_2)$ and sample $20$ values of $\sigma_1$ and $20$ values of $\sigma_2$ uniformly from a base-10 logarithmic distribution from $10^{-1}$ to $10^{1}$, as in \autoref{ss:twomaterials1d}, leading to $20^2$ parameter tuples, and a training set of $280$ $\bs{\sigma}$-tuples, so that, with a batch size of $1$, each epoch consisted of $280$ iterations.
The reference Ritz energy (which acts as the `exact' value) for every $\bs{\sigma}$ was computed with Netgen/NGSolve using piecewise-linear functions on a fine unstructured triangular mesh.
\autoref{fig:2dLshape-training-errors} shows the parametric NN training on the $N=10$ and $N=64$ cases, though we note that the number of `active' elements is variable due to the bottom-right quadrant being removed. For $N=10$, training is carried out for $5\mathord{,}000$ epochs ($1\mathord{,}400\mathord{,}000$ iterations) with a fixed LR of $10^{-2}$. Meanwhile, for $N=64$, we train using Adam with an initial LR of $10^{-2}$ for $50$ epochs ($14\mathord{,}000$ iterations), followed by a LR of $10^{-3}$ for the subsequent $300$ epochs ($84\mathord{,}000$ iterations), and further decreased to $5\cdot10^{-5}$ thereafter, for a total of $12\mathord{,}500$ epochs ($3\mathord{,}500\mathord{,}000$ iterations).

\begin{figure}[htb!]
    \centering
    \begin{subfigure}[t]{0.4\textwidth}
         \resizebox{\textwidth}{!}{%
         \input{figures/example_II_2b/training-error-N10.tex}
         }
        \caption{Case $N=10$.}
        \label{fig:example_II_2b-training-errors-N16}
    \end{subfigure}\mbox{}\hspace{4mm}
    \begin{subfigure}[t]{0.4\textwidth}
         \resizebox{\textwidth}{!}{%
         \input{figures/example_II_2b/training-error-N64.tex}%
         }
        \caption{Case $N=64$.}
        \label{fig:exampleII2b-training-errors-N256}
    \end{subfigure}
    \caption{
    Training histories during $r$-adaptivity for $N=10$ and $N=64$ elements per side in the parametric model problem \eqref{eq:model_problem2} solving $\nabla\cdot(\sigma(\bs{x})\nabla u)=1$ and $u|_{\partial\Omega}=0$ on a multi-material L-shape domain with $\sigma(\bs{x})$ defined by \eqref{eq:materialproLshape}.
    }
    \label{fig:2dLshape-training-errors}
\end{figure}

\autoref{tab:errors_Lshape} shows the error statistics resulting from the $r$-adaptive parametric NN and from uniform meshes.
For the test set with $N=64$, the $r$-adaptive method gives a mean relative error of $0.032$ and a maximum relative error of $0.054$, which is better than the $0.045$ and $0.103$ obtained with uniform grids, respectively.
Similar improvements are seen for the training set and for $N=10$.
Thus, the parametric NN generates well-adapted meshes, enhancing accuracy across different parameters in this challenging two-dimensional scenario.

\begin{table}[htb!]
    \centering
    \caption{Relative error statistics for the uniform and the parametric $r$-adapted solutions associated with  \autoref{fig:2dLshape-training-errors}.}
    \label{tab:errors_Lshape}
    \pgfplotstableread[col sep=comma]{data/example_II_2b/Error_result_N9_15M.csv}\dataNTen
    \pgfplotstableread[col sep=comma]{data/example_II_2b/Error_result_N64.csv}\dataNSixtyFour

    \begin{tabular}{c c c c c c}
        \toprule
        $N$ & Dataset &
        \makecell{$r$-adaptive \\ mean $e_\theta^\sigma$} &
        \makecell{$r$-adaptive \\ max $e_\theta^\sigma$} &
        \makecell{uniform grid \\ mean $e_h^\sigma$} &
        \makecell{uniform grid \\ max $e_h^\sigma$} \\
        \midrule

        \multirow{2}{*}{10}
        & \pgfplotstablegetelem{0}{Set}\of\dataNTen
        \pgfplotsretval &
        \pgfplotstablegetelem{0}{Mean Relative Error}\of\dataNTen
        \pgfmathprintnumber[fixed, zerofill, precision=4]{\pgfplotsretval} &
        \pgfplotstablegetelem{0}{Max Relative Error}\of\dataNTen
        \pgfmathprintnumber[fixed, zerofill, precision=4]{\pgfplotsretval} &
        \pgfplotstablegetelem{0}{Uniform grid Mean Relative error}\of\dataNTen
        \pgfmathprintnumber[fixed, zerofill, precision=4]{\pgfplotsretval} &
        \pgfplotstablegetelem{0}{Uniform grid Max Relative error}\of\dataNTen
        \pgfmathprintnumber[fixed, zerofill, precision=4]{\pgfplotsretval} \\
        & \pgfplotstablegetelem{1}{Set}\of\dataNTen
        \pgfplotsretval &
        \pgfplotstablegetelem{1}{Mean Relative Error}\of\dataNTen
        \pgfmathprintnumber[fixed, zerofill, precision=4]{\pgfplotsretval} &
        \pgfplotstablegetelem{1}{Max Relative Error}\of\dataNTen
        \pgfmathprintnumber[fixed, zerofill, precision=4]{\pgfplotsretval} &
        \pgfplotstablegetelem{1}{Uniform grid Mean Relative error}\of\dataNTen
        \pgfmathprintnumber[fixed, zerofill, precision=4]{\pgfplotsretval} &
        \pgfplotstablegetelem{1}{Uniform grid Max Relative error}\of\dataNTen
        \pgfmathprintnumber[fixed, zerofill, precision=4]{\pgfplotsretval} \\
        \midrule

        \multirow{2}{*}{64}
        & \pgfplotstablegetelem{0}{Set}\of\dataNSixtyFour
        \pgfplotsretval &
        \pgfplotstablegetelem{0}{Mean Relative Error}\of\dataNSixtyFour
        \pgfmathprintnumber[fixed,zerofill, precision=4]{\pgfplotsretval} &
        \pgfplotstablegetelem{0}{Max Relative Error}\of\dataNSixtyFour
        \pgfmathprintnumber[fixed,zerofill, precision=4]{\pgfplotsretval} &
        \pgfplotstablegetelem{0}{Uniform grid Mean Relative error}\of\dataNSixtyFour
        \pgfmathprintnumber[fixed,zerofill, precision=4]{\pgfplotsretval} &
        \pgfplotstablegetelem{0}{Uniform grid Max Relative error}\of\dataNSixtyFour
        \pgfmathprintnumber[fixed, zerofill, precision=4]{\pgfplotsretval} 
        \\
        & \pgfplotstablegetelem{1}{Set}\of\dataNSixtyFour
        \pgfplotsretval &
        \pgfplotstablegetelem{1}{Mean Relative Error}\of\dataNSixtyFour
        \pgfmathprintnumber[fixed,zerofill, precision=4]{\pgfplotsretval} &
        \pgfplotstablegetelem{1}{Max Relative Error}\of\dataNSixtyFour
        \pgfmathprintnumber[fixed,zerofill, precision=4]{\pgfplotsretval} &
        \pgfplotstablegetelem{1}{Uniform grid Mean Relative error}\of\dataNSixtyFour
        \pgfmathprintnumber[fixed, zerofill, precision=4]{\pgfplotsretval} &
        \pgfplotstablegetelem{1}{Uniform grid Max Relative error}\of\dataNSixtyFour
        \pgfmathprintnumber[fixed,zerofill, precision=4]{\pgfplotsretval} \\
        \bottomrule
    \end{tabular}
\end{table}

Lastly, for $N=10$, \autoref{fig:LshapecomparisonParvsNonPar} compares the performance of the parametric NN with that of non-parametric NNs on some samples of the test set. Each non-parametric network is trained for $15\mathord{,}000$ epochs using a fixed learning rate of $10^{-2}$. 
In line with other examples above, and likely for the same reasons, the results coming from the parametric NN yield better results (lower values of Ritz energy) than the non-parametric NN trained at the specific value of $\bs{\sigma}$. 
This reflects what appears to be an added (and somewhat counterintuitive) benefit of having a parametric NN for a class of problems, rather than a NN specialized to a specific problem.

\begin{figure}[htb!]
    \centering
    \begin{subfigure}[t]{0.05\textwidth}
        \begin{tikzpicture}[scale=0.9]
        \node[anchor=west] at (0cm,0cm) {};
            \node[anchor=west] at (-0.5cm,7.0cm) {$u_{\theta}$};
            \node[anchor=west] at (-0.5cm,2.0cm) {$u_{\theta}^{\sigma}$};
        \end{tikzpicture}
     \end{subfigure}\mbox{}\hspace{1mm}
    \begin{subfigure}[t]{0.25\textwidth}
        \resizebox{\textwidth}{!}{
        \input{figures/example_II_2b/test1_r_NN}
        }
        \caption{Test sample $\bs{\sigma} = (2.34, 0.1)$:\\
        $e_\theta=0.2193$ and $e_\theta^{\sigma}=0.2157$.}
    \end{subfigure}\mbox{}\hspace{2mm}
    \begin{subfigure}[t]{0.25\textwidth}
        \resizebox{\textwidth}{!}{
        \input{figures/example_II_2b/test2_r_NN}
        }
        \caption{Test sample $\bs{\sigma} = (1.83, 4.83)$:\\
        $e_\theta=0.1979$ and $e_\theta^{\sigma}=0.1699$.}
    \end{subfigure}\mbox{}\hspace{2mm}
    \begin{subfigure}[t]{0.25\textwidth}
    \resizebox{\textwidth}{!}{
        \input{figures/example_II_2b/test3_r_NN}
        }
        \caption{Test sample $\bs{\sigma} = (0.89, 0.55)$:\\
        $e_\theta=0.2147$ and $e_\theta^{\sigma}=0.1906$.}
    \end{subfigure}
    \caption{
    Finite element solutions using a uniform and $r$-adapted mesh on three test samples for $N=10$. In each case, we show the solution $u_\theta$ coming from a case-by-case non-parametric $r$-adaptivity, and the solution $u^\sigma_\theta$ coming from the parametric $r$-adaptivity produced during optimization of \autoref{fig:2dLshape-training-errors}.
    \vspace{-2mm}
    }
    \label{fig:LshapecomparisonParvsNonPar}
\end{figure}

\section{Conclusions and future work}\label{section5}

This work presents a novel $r$-adaptive FEM tailored to solve parametric PDEs using NNs. 
The core of the proposed approach involves using an NN to generate an optimized, parameter-dependent mesh, which is then used by a standard FEM solver to compute the solution for each parameter instance. 
This strategy combines the robustness and reliability guarantees inherent to the FEM with the capabilities of NNs for adaptive node relocation. 
By concentrating mesh elements in regions with high-variation solutions, the method achieves significant accuracy improvements compared to uniform meshes, particularly for problems exhibiting sharp gradients, boundary layers, or singularities. 
The $r$-adaptive method approaches the convergence rate associated with smooth solutions, even when the solution is singular.
The method was implemented on top of JAX to benefit from its automatic differentiation and accelerated linear algebra via just-in-time compilation, but the linear solve associated with the FEM can be implemented outside of this framework with any desirable solver, thus maximizing performance.

The effectiveness of the NN-based $r$-adaptive FEM approach is demonstrated through numerical experiments on one- and two-dimensional parametric Poisson problems, showcasing its ability to robustly handle varying parameters in both the diffusion coefficient and the source term. The training of the NN relies on minimizing the Ritz energy, which currently restricts the method's applicability to symmetric and coercive problems, such as those arising from elliptic self-adjoint problems. 
Furthermore, the current implementation focuses on first-order piecewise-polynomial elements and tensor-product meshes.

Future research will focus on extending its applicability to non-symmetric or indefinite PDEs, potentially by modifying the loss function. Another important avenue is to extend our implementation to high-order polynomials and irregular geometries, which will require the development of more sophisticated mesh generation strategies that remain compatible with the automatic differentiation required by NN training. Finally, we will investigate the integration of goal-oriented adaptive schemes to further refine the mesh optimization process based on specific quantities of interest.

\section*{Acknowledgments}
AJO, CU and DP declare this work has received funding from the following Research Projects/Grants: European Union’s Horizon Europe research and innovation programme under the Marie Sklodowska-Curie Action MSCA-DN-101119556 (IN-DEEP); PID2023-146678OB-I00 funded by MICIU/AEI/10.13039/501100011033 and by FEDER, EU; PID2023-146668OA-I00 funded by MICIU/AEI/10.13039/501100011033 and by FEDER, EU; BCAM Severo Ochoa accreditation of excellence CEX2021-001142-S funded by MICIU/AEI/10.13039/501100011033; Basque Government through the BERC 2022-2025 program; BEREZ-IA (KK-2023/00012) and RUL-ET(KK-2024/00086), funded by the Basque Government through ELKARTEK; Consolidated Research Group MATHMODE (IT1456-22) given by the Department of Education of the Basque Government; BCAM-IKUR-UPV/EHU, funded by the Basque Government IKUR Strategy and by the European Union NextGenerationEU/PRTR. 
DA, VI, IT, FF and MS were partially supported by the National Center for Artificial Intelligence CENIA FB210017, Basal ANID based in Chile. DA, VI and FF were partially supported by the Office of Naval Research (ONR) award N629092312098, and MS was partially funded by FONDECYT Regular (ANID), grant N.~1221189. 
Lastly, all authors acknowledge support by grant N.~FOVI240119 (Fomento a la Vinculaci\'on Internacional para instituciones de Investigaci\'on), ANID Chile.

\bibliographystyle{unsrt}

\end{document}